\documentclass[3p]{elsarticle}
\usepackage{color}
\usepackage{amssymb}
\usepackage[figuresright]{rotating}

\def\R{{\rm I\kern-.17em R}}
\def\Z{{\rm Z\kern-.32em Z}}
\def\ds{\displaystyle}
\usepackage{epsf} 
\usepackage{epsfig}
\usepackage{graphicx}
\usepackage{comment}
\usepackage{amssymb}
\usepackage{mathtools}
\usepackage{algorithm}
\def\ds{\displaystyle}

\def\N{{\rm I\kern-.20em N}}
\def\R{{\rm I\kern-.17em R}}
\def\Z{{\rm Z\kern-.32em Z}}

\def\ds{\displaystyle}
\def\bkC{{\rm \kern.24em\vrule 
width.05em height1.4ex depth-.05ex\kern-.26emC}}
\DeclareMathOperator*{\argmin}{argmin}
\begin{document}

\begin{frontmatter}

\title{Optical flow with fractional order regularization: variational model  and solution method}

\author[cmuc]{Somayeh Gh. Bardeji\fnref{fn1}}
\author[cmuc]{Isabel N. Figueiredo\fnref{fn1}}
\author[cmuc]{Erc\'\i lia Sousa\corref{cor1}\fnref{fn1}}
\cortext[cor1]{Corresponding author}

\address[cmuc]{CMUC, Department of Mathematics, University of Coimbra, 
3001-501 Coimbra, Portugal}
\fntext[fn1]{Research supported in part  by
the Portuguese National Funding Agency for Science, Research and Technology (FCT) under Project
PTDC/MATNAN/0593/2012 and by CMUC -- UID/MAT/00324/2013, funded by the Portuguese Government through 
FCT/MEC and co-funded by the European Regional Development Fund through the Partnership Agreement PT2020.}
\ead{ecs@mat.uc.pt}

\begin{abstract}

An optical flow  variational model is proposed for a sequence of images
defined on a domain in $\mathbb{R}^2$.  We introduce a regularization term given by  the $L^1$ 
norm of a fractional differential  operator. 
To solve the minimization problem 
we apply the split Bregman method.
Extensive experimental results, with performance evaluation,
are presented to demonstrate the effectiveness of the new model and method
and to show that our algorithm performs favorably in comparison to another 
existing method.  We also discuss the influence of the 
order $\alpha$ of the fractional operator in the estimation of the optical flow, 
for $0 \leq \alpha \leq 2$. We observe that the values of $\alpha$ for which the method 
performs better depends on the geometry and texture complexity of the image.
  Some extensions of our algorithm are also discussed.
\end{abstract}

\begin{keyword}
Optical flow, variational models, fractional derivatives, finite
differences.
\end{keyword}
\end{frontmatter}



\section{Introduction}
Optical flow is a tool for detecting and analyzing motion in a
sequence of images. The underlying idea is to depict the displacement of patterns in the
image sequence as a vector field, named the optical flow vector field, generating the corresponding
displacement function.
In their seminal paper, Horn and Schunck \cite{hor} suggested a variational method
 for the computation of the optical flow vector field. In this approach the goal is to minimise
 an energy functional consisting of a similarity term (or data term) and a regularity term:
 $$
\argmin_{u \in \cal{H}} E({\bf u}) = 
\argmin_{{\bf u} \in \cal{H}}({\cal{R}}({\bf u})+{\cal{S}}({\bf u})).
$$
The space $\cal{H}$ denotes an admissible space of vector fields, $\cal{R}$ denotes the regularity term for the
vector field ${\bf u}$, and  $\cal{S}$ denotes the similarity term that
depends on the data image sequence. In particular the functional is of the form \cite{hor}
\begin{equation}
E(\mathbf{u})= \beta^2\int_\Omega  (|\nabla u_1|^2 + |\nabla u_2|^2) d\Omega          
+ \int_\Omega (I_1(\mathbf{x}+\mathbf{u}(\mathbf{x}))- I_0(\mathbf{x}))^2d\Omega.
\label{hs}
\end{equation}
Here, $I_0$ and $I_1$ is the image pair, ${\bf u} = (u_1({\bf x}), u_2({\bf
  x}))^T$ is the two-dimensional displacement
field and $\beta$ is a fixed parameter. The first term (regularization term)
penalizes high variations in ${\bf u}$ to obtain smooth displacement fields. The second
term (data term) is also known as the optical flow constraint. It assumes, that
the intensity values of $I_0({\bf x})$ do not change during its motion to 
$I_1({\bf x} + {\bf u}({\bf x}))$.
Horn and Schunck \cite{hor}  observed that $\beta^2$ plays a significant role only for areas where
the brightness gradient is small, preventing haphazard adjustments to the
estimated flow velocity.  Disadvantages of this model consist of 
not preserving discontinuities in the flow field and of not handling outliers efficiently.
To  overcome  the difficulties presented by the Horn-Schunck functional, 
several extensions and improvements have been developed  \cite{wei}. 

In \cite{zac2007} the optical flow model proposed
consists in considering an $L^1$ norm in the regularizing term and the 
similarity  term is substantially changed by introducing an auxiliary variable
${\bf v}$.  The process is a result of first changing the
 quadratic factors that appeared in the classical method (\ref{hs}), obtaining
an energy functional which is the sum of the total variation of ${\bf u}$ and an $L^1$ term:
\begin{equation}
E({\bf u}) = \int_\Omega |\nabla {\bf u}| d\Omega +
\lambda\int_\Omega|\rho(\mathbf{u})|d\Omega,
\label{Ei}
\end{equation}
where $|\nabla {\bf u}|=|\nabla u_1|+|\nabla u_2|$ and 
the image residual denoted by $\rho(\mathbf{u})$
(we omitted the explicit dependency on $\mathbf{u}^0$ and $\mathbf{x}$)
is given by
\begin{equation}				
\rho(\mathbf{u})=\nabla I_1(\mathbf{x} +
\mathbf{u^0}).(\mathbf{u}-\mathbf{u}^0) + I_1(\mathbf{x}+             
\mathbf{u}^0) - I_0(\mathbf{x}).
\label{rho}
\end{equation}
The vector ${\bf u}^0$ is a given disparity map and
the functional was obtained for a fixed ${\bf u}^0$ and using the linear
approximation for $I_1({\bf x}+{\bf u})$ near ${\bf x}+{\bf u}^0$.

Secondly, a convex relaxation term is introduced \cite{zac2007}  in order to minimize this energy
functional 
efficiently obtaining  
\begin{equation}
E_\theta(\mathbf{u},\mathbf{v})= \int_\Omega \left\{|\nabla {\bf u}|
+ \frac{1}{2\theta} |\mathbf{u}-\mathbf{v}|^2
+ \lambda|\rho(\mathbf{v})|\right\} d\Omega,
\label{Etheta}
\end{equation}		
where
$\theta$ is a small constant, such that ${\bf v}$ is a close
approximation of ${\bf u}$. Setting $\theta$ very small forces the minimum of
$E_\theta$
to occur when ${\bf u}$ and ${\bf v}$ are nearly equal, reducing the energy
(\ref{Etheta})
to the original energy 
(\ref{Ei}).

Many approaches for optical flow computation replace the nonlinearity
intensity profile $I_1({\bf x} + {\bf u})$ by a first Taylor approximation
to linearize the problem locally as in the case presented above.
Since such approximations are only
valid for small motions, in the presence of large
displacements, the method fails when the gradient of the image is not smooth
enough.  This means that additional techniques are required to determine the optical
flow correctly. Therefore  an iterative warping is applied in the implementation to
compensate for image nonlinearities. A multiscale strategy is also included
to allow disparities between the images.


In this work we propose an optical flow  model   for a sequence of
images defined on a domain in $\mathbb{R}^2$ which consists of a modification
of the  model introduced in \cite{zac2007}, by considering for the 
  regularization term  the $L^1$ norm of a fractional derivative operator \cite{che2013}.
The numerical method developed to solve the minimization problem involves
a multiscale strategy \cite{mei2013}
and the split Bregman method described in \cite{gol2009}. 
The effectiveness of the new model and  numerical approach  is shown 
by presenting  experimental results that
use the test sequences  available in the Middlebury benchmark database designed by \cite{bak}.
We also compare its performance with other  existing numerical method.

In the next section we present the variational method and in Section 3 we describe  
the numerical approach which includes the split Bregman method,
Euler Lagrange equations, a shrinkage operator, a thresholding operator and finite differences.
In Section 4 several experiments are shown
and we end with  some conclusions and general comments in Section 5.

\section{Problem formulation}

We propose a generalised method that involves fractional derivatives in		
the regularisation term.
Recently fractional derivatives have been brouhgt to the field of 
image processing and fractional differentiation based methods
have been demonstrating advantages over already existing methods, 
see for instance \cite{cha2013, che2013, rom2008, zha2012}.
We start to introduce the definition of fractional		
derivative. 

The left Riemann-Liouville derivative of order $\alpha$, for a scalar function $u$,
is defined by		
\begin{equation}	
D_{-}^{\alpha} u(t) =\frac{1}{\Gamma(m-\alpha)} \frac{d^m}{dt^m}		
\int_{a}^{t} u(\tau) (t-\tau)^{m-\alpha-1}d\tau, \quad m-1<\alpha<m,
\label{frac}		
\end{equation}
for $a \leq t \leq b$,
where $m$ is a positive integer and $\Gamma$ denotes the Gamma function
$$
\ds{\Gamma(z)=\int_{0}^{\infty} x^{z-1} {\rm e}^{-x}dx},
$$
that satisfies the property $\Gamma(z+1)=z\Gamma(z)$.
In this work we are interested in  $0\leq \alpha\leq 2$.  
When $\alpha$ goes to $m$  the operator $D_{-}^{\alpha} u(t)$ becomes
the integer order derivative $u^{(m)}(t)$. Hence, the fractional
derivative  is seen as a generalization of the classical derivative.

The motivation to include the fractional operator has to do with  its capability
to change continuously the regularization operator depending on the choice of
the value of $\alpha$. In particular, for $\alpha=0,1,2$ it represents  the function, 
the first order derivative and the second order derivative respectively.
This allows to choose the most suitable $\alpha$ (regularization operator)
for different types of images according to its  high or low texture,
the presence of motion discontinuities, flat, corners or edges.

Let $u$ be a  scalar function defined in $[a,b] \times [c,d]$. 
For $(x,y) \in (a,b)\times (c,d)$, we define	
\begin{equation}
D_{x-}^{\alpha} u(x,y) =\frac{1}{\Gamma(m-\alpha)} \frac{d^m}{dx^m}		
\int_{a}^{x} u(\tau,y) (x-\tau)^{m-\alpha-1}d\tau
\label{dxn}	
\end{equation}	
\begin{equation}
D_{y-}^{\alpha} u(x,y) =\frac{1}{\Gamma(m-\alpha)} \frac{d^m}{dy^m}	
\int_{c}^{y} u(x,\tau) (y-\tau)^{m-\alpha-1}d\tau.		
\label{dyn}	
\end{equation}	
The left
Riemann-Liouville fractional operator $\nabla_{-}^\alpha u$   denotes
$\nabla_-^\alpha u= (D^\alpha_{x-}u,D^\alpha_{y-} u)$ with euclidean norm
$|\nabla_-^\alpha u| = \sqrt{ (D^\alpha_{x-} u)^2 +(D^\alpha_{y-} u)^2}.
$

A generalised form of  the energy (\ref{Etheta}) is 
\begin{equation}	
E^\alpha_\theta(\mathbf{u},\mathbf{v})= 
\int_\Omega \left\{|\nabla_{-}^{\alpha} {\bf u}|
+ \frac{1}{2\theta} |\mathbf{u}-\mathbf{v}|^2+ \lambda		
|\rho(\mathbf{v})|\right\}d\Omega,
\label{Ethetaalpha}		
\end{equation}
where ${\bf u}({\bf x})=(u_1({\bf x}),u_2({\bf x}))$ and
$
|\nabla_{-}^{\alpha} {\bf u}| = |\nabla_-^\alpha u_1| + |\nabla_-^\alpha u_2|.
$	
The minimisation of the energy
$E^\alpha_\theta$ can be performed by  alternating steps \cite{zac2007}
and updating either $\mathbf{u}$ or $\mathbf{v}$ at each iteration, that is,
first we fix $\mathbf{v}$, and solve 			
\begin{eqnarray}			
\min_{\bf u} \int_\Omega \left\{|\nabla _{-}^\alpha {\bf u}|
+ \frac{1}{2\theta}|\mathbf{u}-\mathbf{v}|^2\right\}d\Omega,
\label{first_minimization}		
\end{eqnarray}
and secondly we fix $\mathbf{u}$ and solve			
\begin{eqnarray}\label{second_minimization}		
\min_{\mathbf{v}}\int_\Omega \left\{\frac{1}{2\theta}|\mathbf{u}-\mathbf{v}|^2 
+ \lambda|\rho(\mathbf{v})|\right\}d\Omega.
\label{second_minimization}		
\end{eqnarray}
If {\bf v} is fixed, the functional is convex in {\bf u}. Therefore,
a global minimizer {\bf u} can be computed efficiently.
If {\bf u} is
fixed, the functional has only a pointwise dependency
on {\bf v}. Therefore, it can be minimized globally
with respect to {\bf v} by a complete search.

Note that when $\alpha=1$ we get 
$\nabla _-^1{\bf u}=\nabla {\bf u}$ and in this case problem (\ref{Ethetaalpha})
reverts to problem (\ref{Etheta}).
A way to solve (\ref{first_minimization}) was proposed in \cite{cha2004}
(for the case $\alpha=1$) which uses a dual formulation of (\ref{first_minimization}) to derive an
efficient and globally convergent scheme. In our work to solve this equation
we use the split Bregman tecnhique \cite{gol2009}.
Equation (\ref{second_minimization}) is solved using the approach presented in
\cite{san2013, zac2007} and in Section 3.3 we report briefly this known tecnhique.

\section{Numerical method}

In this section we describe our proposal on how to estimate the optical flow by solving the two minimisation problems
presented  in the previous section.
It consists essentially in the application of the split Bregman technique, the derivation of 
Euler Lagrange equations and the use of finite difference approximations.
The first minimisation problem is discussed, in section 3.1,
for the  particular case $\alpha=1$ and then  the general case, for
$0 \leq \alpha\leq 2$, is described in section 3.2.  
The second minimization problem does not
depend on $\alpha$ and therefore is solved equally for all cases
and presented in section 3.3. In section 3.4 we describe in detail the implementation of
the algorithm.

\subsection{Solving problem (\ref{first_minimization}) with split Bregman method for $\alpha=1$}
Consider the minimization problem (\ref{first_minimization}) when
$\alpha=1$. To find its solution  we can  solve
the following problem:
for any fixed $(v_1, v_2)$, search for $u_\ell$, $\ell=1,2$
\begin{equation}			
\min_{u_\ell} \int_\Omega \left\{|\nabla u_\ell| + \frac{1}{2\theta}|u_\ell -
  v_\ell|^2
\right\}d\Omega,
\quad \ell=1,2.
\label{u1u2}			
\end{equation}
We propose to solve this problem with  split Bregman method, since
 (\ref{u1u2}) belongs to the general class of problems discussed
in \cite{gol2009},
$
\ds{\min_{u} \{\|\phi(u)\|_{L^1(\Omega)}  + H(u)\}},
$
where $\|\cdot \|_{L^1(\Omega)}$ denotes the $L^1$ norm and both
$\|\phi(u)\|_{L^1(\Omega)}$
and $H(u)$  are convex functions. A brief explanation of the split
Bregman method is given in what follows. 

We first replace problem (\ref{u1u2}) 
by the constrained optimization problem 
\begin{equation}		
\min_{u_\ell}\left\{ \|d_\ell\|_{L^1(\Omega)}  + \frac{1}{2\theta}\|u_\ell-v_\ell\|_{L^2(\Omega)}^2\right\}, \quad \ell=1,2
\label{conu1u2}	
\end{equation}
subject to $d_\ell=\nabla u_\ell$, $\ell=1,2$. Then to get an unconstrained
problem, a $L^2$ penalty term is added
\begin{equation}
\min_{d_\ell,u_\ell} \left\{\|d_\ell\|_{L^1(\Omega)}  + \frac{1}{2\theta}\|u_\ell - v_\ell\|_{L^2(\Omega)}
+ \frac{\lambda_\ell}{2}\|d_\ell - \nabla u_\ell\|_{L^2(\Omega)}^2 \right\}.
\label{unconu1u2}		
\end{equation}
The problem is then modified to get exact enforcement of the constraint using
a Bregman iteration \cite{ber1967}. This leads to the split Bregman method
that consists of solving the following problem. 
For $k=1,2, \dots$; $ \ \ \ell=1,2$
\begin{equation}
\left\{\begin{array}{l}
(u_\ell^{k+1}, d_\ell^{k+1}) = \displaystyle{\min_{d_\ell, u_\ell}} \left\{\|d_\ell\|_{L^1(\Omega)} +
\frac{1}{2\theta}\|u_\ell - v_\ell\|_{L^2(\Omega)}^2 + \frac{\lambda_\ell}{2}\|d_\ell -
\nabla u_\ell - b_\ell^k\|_{L^2(\Omega)}^2 \right\}\\
\\
b_\ell^{k+1} = b_\ell^k - d_\ell^{k+1}  + \nabla u_\ell^{k+1}.
\end{array}
\right.
\label{csb}
\end{equation}
Problem (\ref{csb}) is solved by alternate iterative minimization, meaning that two
steps are performed. First for a fixed $d_\ell$ the minimization is done with
respect to $u_\ell$, that is, 
\begin{equation}
u_\ell^{k+1} = \min_{u_\ell} \left\{\frac{1}{2\theta} \|u_\ell - v_\ell\|_{L^2(\Omega)}^2 
+ \frac{\lambda_\ell}{2}\|d_\ell^k-\nabla u_\ell - b_\ell^k\|_{L^2(\Omega)}^2\right\}.
\label{step1}
\end{equation}
Secondly, for a fixed $u_\ell$ the minimization is done with respect to $d_\ell$,
that is,
\begin{equation}
d_\ell^{k+1}  = \min_{d_\ell} \left\{\|d_\ell\|_{L^1(\Omega)}
+ \frac{\lambda_\ell}{2}\|d_\ell-\nabla u_\ell^{k+1} - b_\ell^k\|_{L^2(\Omega)}^2\right\}.
\label{step2}
\end{equation}

To solve (\ref{step1}) we derive the optimality condition for $u_\ell^{k+1}$,
that consists on the following  Euler-Lagrange equations. For
$\ell=1,2$,
\begin{equation}
\left(\frac{1}{\theta}  - \lambda_\ell \Delta\right) u_\ell^{k+1} 
= \frac{1}{\theta}  v_\ell - \lambda_\ell \ {\rm div} (d_\ell^k - b_\ell^k) \qquad \mbox{in} \ \Omega
\label{eulerlagrange}
\end{equation}
and
subject to the natural boundary conditions 
\begin{equation}
\frac{\partial u_{\ell}}{\partial \eta}  = (d_{\ell}^k - b_{\ell}^k )\cdot \eta \qquad \mbox{on} \ \partial \Omega.
\label{bceulerlagrange}
\end{equation}
Here,  $\Delta$ and ${\rm div}$ denote the Laplace and divergent operators
respectively, ``$\cdot$'' denotes the inner product in $\R^2$
and $\eta$ is the unit outward normal to $\partial \Omega$.

We use finite differences to approximate  the derivatives
and since the resulting system is diagonally dominant, it
can be solved efficiently with the Gauss-Seidel iterative method.
The solution $ (u_\ell)_{i,j}^{k+1}$  at each pixel $(i,j)$ in $\Omega$ is given by
\begin{equation}		
(u_\ell)_{i,j}^{k+1} =\frac{1}{\displaystyle{{1}/{\theta}} + 4\lambda_{\ell}}
\left[\lambda_{\ell}(U_\ell)^{k+1}_{i,j}  
+ \frac{1}{\theta} (v_\ell)_{i,j} -   \lambda_{\ell}  \ \mbox{div}( d^k  - b^k)_{i,j}  \right],
\end{equation}
where $(U_\ell)^{k+1}_{i,j} $ denotes
$	
(U_\ell)^{k+1}_{i,j}:= (u_\ell)_{i-1,j}^{k}+ (u_\ell)_{i+1,j}^{k+1} +
(u_\ell)_{i,j-1}^{k} + (u_\ell)_{i,j+1}^{k+1}
$
and to approximate the divergence operator backward finite differences have been applied leading to
$$
\mbox{div}( d^k  - b^k)_{i,j} = 2d_{i,j}^k-d_{i-1,j}^k-d_{i,j-1}^k-2b_{i,j}^k+b_{i-1,j}^k+b_{i,j-1}^k.
$$

To solve (\ref{step2}) a shrinkage operation can be used at each point
$(i,j)$ such that
\begin{equation}		
d_\ell^{k+1} = shrink\left(\nabla u_\ell^{k+1} + b_\ell^k, \frac{1}{\lambda_\ell}\right),
\label{shrink}		
\end{equation}
where 
$$
\quad \ds{shrink(x,\gamma)= \frac{x}{|x|} \ \max(|x|-\gamma, 0)}, \ z, \gamma \in \R.
$$

\subsection{Solving problem (\ref{first_minimization}) with split Bregman method for $0\leq
  \alpha\leq 2$}

In this section we solve the minimization problem (\ref{first_minimization}) for
$0\leq \alpha\leq2$, that is, 
for any fixed $(v_1, v_2)$, we search for the minimizer $(u_1,u_2)$ of the problem
\begin{equation}			
\min_{u_\ell} \int_\Omega \left\{|\nabla_{-}^\alpha u_\ell|
+ \frac{1}{2\theta}|u_\ell - v_\ell|^2\right\}d\Omega,
\quad \ell=1,2.
\label{u1u2alpha}			
\end{equation}
Similarly to what we have done in the previous section,
we first replace (\ref{u1u2alpha}) 
by the constrained optimization problem 
\begin{equation}		
\min_{u_\ell} \left\{\|d_\ell\|_{L^1(\Omega)}  + \frac{1}{2\theta}\|u_\ell-v_\ell\|_{L^2(\Omega)}^2\right\}, \quad \ell=1,2
\label{conu1u2alpha}	
\end{equation}
now subject to $d_\ell=\nabla_{-}^{\alpha} u_\ell$, $\ell=1,2$. Then to get an unconstrained
problem a $L^2$ penalty term is added
\begin{equation}
\min_{d_\ell,u_\ell} \left\{\|d_\ell\|_{L^1(\Omega)}  + \frac{1}{2\theta}\|u_\ell - v_\ell\|_{L^2(\Omega)}
+ \frac{\lambda_\ell}{2}\|d_\ell - \nabla_{-}^{\alpha}  u_\ell\|_{L^2(\Omega)}^2 \right\}.
\label{unconu1u2alpha}		
\end{equation}
The problem is then modified to get exact enforcement of the constraint using
a Bregman iteration \cite{ber1967}. This leads to the split Bregman method
that consists of solving the following problem. 
For $k=1,2, \dots$; $ \ \ \ell=1,2$
\begin{equation}
\left\{\begin{array}{l}
(u_\ell^{k+1}, d_\ell^{k+1}) = \displaystyle{\min_{d_\ell, u_\ell}} \left\{\|d_\ell\|_{L^1(\Omega)} +
\frac{1}{2\theta}\|u_\ell - v_\ell\|_{L^2(\Omega)}^2 + \frac{\lambda_\ell}{2}\|d_\ell -
\nabla_{-}^{\alpha}  u_\ell - b_\ell^k\|_{L^2(\Omega)}^2 \right\}\\
\\
b_\ell^{k+1} = b_\ell^k - d_\ell^{k+1}  + \nabla_{-}^{\alpha}  u_\ell^{k+1}.
\end{array}
\right.
\label{csbalpha}
\end{equation}
A solution to problem (\ref{csbalpha}) can be obtained by alternate iterative minimization.
First for a fixed $d_\ell$ the minimization is done with
respect to $u_\ell$, that is, 
\begin{equation}
u_\ell^{k+1} = \min_{u_\ell} \left\{\frac{1}{2\theta} \|u_\ell - v_\ell\|_{L^2(\Omega)}^2 
+ \frac{\lambda_\ell}{2}\|d_\ell^k-\nabla_{-}^{\alpha} u_\ell - b_\ell^k\|_{L^2(\Omega)}^2\right\}.
\label{step1alpha}
\end{equation}
Secondly for a fixed $u_\ell$ the minimization is done with respect to $d_\ell$,
that is,
\begin{equation}
d_\ell^{k+1}  = \min_{d_\ell} \left\{\|d_\ell\|_{L^1(\Omega)}
+ \frac{\lambda_\ell}{2}\|d_\ell-\nabla_{-}^{\alpha} u_\ell^{k+1} - b_\ell^k\|_{L^2(\Omega)}^2\right\}.
\label{step2alpha}
\end{equation}
For this problem, the rectangular domain $[a,b]\times [c,d]$  is an 
extension of the image pixels domain, obtained by padding around the image.
Dirichlet boundary conditions, for $u_\ell$, $\ell=1,2$,
are imposed on the padding region.

We introduce now		
the definition of right fractional derivative and the property of integration
by parts, since they will be needed in what follows.
The right Riemann-Liouville  derivative or order $\alpha$, 
for a scalar function $u$, for $a \leq t \leq b$, is defined by			
$$			
D_{+}^{\alpha} u(t) =\frac{(-1)^m}{\Gamma(m-\alpha)} \frac{d^m}{dt^m}			
\int_{t}^{b} u(\tau) (\tau-t)^{m-\alpha-1}d\tau, \quad  m-1\leq \alpha \leq m.
$$
The spaces of functions 
$I_{a^+}^\alpha(L_p)$ and $I_{b^-}^\alpha(L_p)$  for $\alpha >0$ and
$1 \leq p < \infty$ are defined by
\begin{eqnarray}
I_{a^+}^\alpha(L_p) & = & \left \{ f: f(x)=\frac{1}{\Gamma(\alpha)}\int_{a}^{x} 
\frac{\phi(t)}{(x-t)^{1-\alpha}} dt, \ x >a, \ \phi \in L_p(a,b) \right\}
\\
I_{b^-}^\alpha(L_p)& = & \left \{ f: f(x)=\frac{1}{\Gamma(\alpha)}\int_{x}^{b} 
\frac{\psi(t)}{(t-x)^{1-\alpha}} dt, \ x < b, \ \psi \in L_p(a,b) \right\},
\end{eqnarray}
where $L_p$ denotes the space of $p$ integrable functions in $(a,b)$.
We have the following integration by parts result.

\newtheorem{thm}{Theorem}			
\begin{thm} (\cite[page 76]{kil2006},
\cite[pages 34,46]{sam1993}):
Let $\alpha >0$, $p\geq 1$, $q \geq 1$, and $1/p+1/q\leq 1+\alpha$			
($p \neq 1, q \neq 1$ in the case when $1/p+1/q=1+\alpha$).			
If $f \in I_{b-}^\alpha (L_p)$ and			
$g \in I_{a+}^\alpha (L_q)$			
then			
$$			
\int_{a}^{b} f(x)(D_{-}^\alpha g)(x) dx =\int_{a}^{b} g(x)(D_{+}^\alpha f)(x) dx.			
$$
\end{thm}

See Theorem 2.3 of \cite{sam1993} on page 43, for necessary 
and sufficient conditions	
for $f \in I_{a^+}^\alpha(L_1)$ and Theorem 2.4 on page 45, 
to see what happens if		
$f \notin I_{a^+}^\alpha(L_1)$.

Let $u$ be a  scalar function defined in $[a,b] \times [c,d]$. 
For $(x,y) \in (a,b)\times (c,d)$, we define	
\begin{equation}		
D_{x+}^{\alpha} u(x,y) =\frac{(-1)^m}{\Gamma(m-\alpha)} \frac{d^m}{dx^m}		
\int_{x}^{b} u(\tau,y) (\tau-x)^{m-\alpha-1}d\tau	
\label{dxm}	
\end{equation}	
\begin{equation}
D_{y+}^{\alpha} u(x,y) =\frac{(-1)^m}{\Gamma(m-\alpha)} \frac{d^m}{dy^m}	
\int_{y}^{d} u(x,\tau) (\tau-y)^{m-\alpha-1}d\tau.		
\label{dym}	
\end{equation}	
The fractional derivatives $D_{x-}^{\alpha} u$ and $D_{y-}^{\alpha} u$
have been defined previously in (\ref{dxn}) and (\ref{dyn}) respectively.

Using the integration by parts property, the optimality condition for
$u_\ell$, $\ell=1,2$ is given by the differential equation
\begin{equation}			
\frac{1}{\theta}u_\ell  + \lambda_\ell\left[			
D^{\alpha}_{x+}D^\alpha_{x-} u_{\ell} + D^{\alpha}_{y+}D^\alpha_{y-} u_{\ell}
\right] 
=\frac{1}{\theta} v_\ell  + \lambda_\ell\left[D^{\alpha}_{x+}(d_{1\ell}^k -
  b_{1\ell}^k) 
+ D^{\alpha}_{y+}(d_{2\ell}^k  - b_{2\ell}^k)\right].
\label{elfrac}			
\end{equation}
A discussion on the solutions of the fractional Euler-Lagrange equations can be seen in
\cite{agr2002} or for a more recent work, we refer also, for instance, to \cite{bou2015}.

To approximate the fractional derivatives we use a standard discretization
\cite{pod1999}. Let $(i,j)$ denote an arbitrary pixel in the image domain.
Then, for 
$i=0,1, \dots, N_x$ and $j=0,1,\dots, N_y$,
\begin{eqnarray}		
D^\alpha_{x-} u(i,j)\approx \sum_{k=0}^{i} w_k^{(\alpha)} u(i-k,j), &&	
D^{\alpha}_{x+} u(i,j)\approx\sum_{k=0}^{N_x-i}w_k^{(\alpha)} u(i+k,j)\\
D^\alpha_{y-} u(i,j)\approx\sum_{k=0}^j w_k^{(\alpha)} u(i,j-k), &&	
D^{\alpha}_{y+} u(i,j)\approx\sum_{k=0}^{N_y-j}w_k^{(\alpha)} u(i,j+k),
\end{eqnarray}
where
the coefficients $w_k^{(\alpha)}$ can be obtained by the recurrence formula
$$
w_0^{(\alpha)}=1, \quad w_k^{(\alpha)}
=\left(1-\frac{\alpha+1}{k}\right)w_{k-1}^{(\alpha)},
\quad  k=1, 2, \cdots,N.
$$
A recent discussion on the order of accuracy of the previous
discretizations can be seen
for instance in \cite{sou2012}.

An approximation for the composition of the operators that appear in (\ref{elfrac})
is 
\begin{eqnarray}		
D^{\alpha}_{x+}D^\alpha_{x-} u(i,j)\approx
\sum_{k=i}^{N_x} w_{k-i}^{(\alpha)}\sum_{p=0}^k w_p^{(\alpha)}  u(k-p,j), &&
D^{\alpha}_{y+}D^\alpha_{y-} u(i,j)\approx
\sum_{k=j}^{N_y} w_{k-j}^{(\alpha)}\sum_{p=0}^k w_p^{(\alpha)}  u(i,k-p).		
\end{eqnarray}
Hence, a discrete version of (\ref{elfrac}) can be represented by		
\begin{eqnarray}		
(u_\ell)_{i,j}^{k+1}&=&\frac{\lambda_{\ell}}{a_\ell}
\left[\sum_{p=i}^{N_x}w^{(\alpha)}_{p-i}\sum_{\ell=0,\ell\neq i}^p w^{(\alpha)}_\ell u^k_{\ell-i,j}
+\sum_{p=j}^{N_y} w^{(\alpha)}_{p-j}\sum_{\ell=0,\ell\neq j}^p w^{(\alpha)}_\ell
u^k_{i,\ell-j}
 +\frac{1}{\theta\lambda_{\ell}} (v_\ell)_{i,j} \right.\\\nonumber		
&& \left. +\sum_{p=0}^{N_x-i}w^{(\alpha)}_{p}(d_{1}^k - b_{1}^k)_{i+p,j} 
+ \sum_{p=0}^{N_y-j}w^{(\alpha)}_{p}(d_{2}^k  - b_{2}^k)_{i,j+p}\right],		
\end{eqnarray}	
for
$$ \ a_\ell=\ds{\frac{1}{\theta} +\lambda_\ell\left(
\sum_{p=i}^{N_x-1}(w^{(\alpha)}_{p-i})^2
 +\sum_{p=j}^{N_y-1}(w^{(\alpha)}_{p-j})^2\right).}$$

To solve (\ref{step2alpha})
a shrinkage operation can be used at each point
$(i,j)$, similarly to what has been done in (\ref{shrink}), that is,
\begin{equation}		
d_\ell^{k+1} = shrink\left(\nabla_{-}^\alpha u_\ell^{k+1} + b_\ell^k, \frac{1}{\lambda_\ell}\right).
\label{shrinkalpha}		
\end{equation}

\subsection{Solving problem (\ref{second_minimization})}

The second minimization problem  (\ref{second_minimization}) is solved by the
same method presented in \cite{san2013, zac2007}, where the solution 
 is given by the
following thresholding step:
\begin{equation}		
\mathbf{v}^{k+1} = \mathbf{u}^{k+1} + TH(\mathbf{u}^{k+1},\mathbf{u}^0),
\label{iterv}		
\end{equation}
with the thresholding operator
\begin{equation}	
{\rm TH}(\mathbf{u},\mathbf{u}^0) = 
\left\{\begin{array}{ll}		
\lambda \theta \nabla I_1(\mathbf{x} + \mathbf{u}^0) 
&{\rm if} \ \rho(\mathbf{u}) < -\lambda \theta |\nabla
I_1(\mathbf{x} + \mathbf{u}^0)|^2\\		
-\lambda \theta \nabla I_1(\mathbf{x} + \mathbf{u}^0)         
& {\rm  if} \ \rho(\mathbf{u}) >  \lambda \theta |\nabla
I_1(\mathbf{x} + \mathbf{u}^0)|^2\\		
-\rho(\mathbf{u}) {\displaystyle \frac{\nabla I_1(\mathbf{x} 
+ \mathbf{u}^0)}{|\nabla I_1(\mathbf{x} + \mathbf{u}^0)|^2}}
& {\rm if} \ \rho(\mathbf{u}) \leq  
\lambda \theta |\nabla I_1(\mathbf{x} + \mathbf{u}^0)|^2.
\end{array}		
\right.
\label{thresholding}		
\end{equation}

The input of the algorithm is a pair of images $I_0({\bf x})$ and $I_1({\bf
  x})$  with ${\bf x}=(i,j)$ the pixel index. The output is a vector field
${\bf u}({\bf x})=(u_1({\bf x}),u_2({\bf x}))$. The residual 
$\rho(\mathbf{u})$ is a scalar field, that is, a gray valued
image. The vector field ${\bf u}^0$ must be close to ${\bf u}$. It is given by
the enclosing multiscale procedure and it is zero at the coarsest level.
The gradient of the image $\nabla I_1$ is approximated with central differences along each
direction and at the borders Neumann boundary conditions are assumed.
To warp the image $I_1$ by a flow field ${\bf u}^0$, we evaluate 
$I_1({\bf x}+{\bf u}^0({\bf x}))$ using bicubic interpolation.
A more detailed information about the implementation of the algorithm is given
in the next section.

\subsection{Implementation of the proposed algorithm}

A variety of approaches have been used to improve the convergence rate
of optical flow algorithms.  We apply a coarse to fine strategy, common to many optical flow
algorithms, that consists of building image pyramids \cite{bla1996,enk1986,san2013}.
The optical flow is first computed on the top level (fewest pixels) and then
upsampled and used to initialize the estimate at the next level. Computation
at the higher levels in the pyramid involves fewer unknows and therefore is
faster. The initialization at each level from the previous level also means
that fewer iterations are required at each level to reach a certain accuracy.
Incremental warping of the flow between pyramid levels helps to keep the flow
update at any given level small (fewer pixels) and when combined
with incremental warping and updating within a level, the method 
becomes very effective.

To create the pyramid of images we follow a standard strategy \cite{mei2013}.
The pyramid is built by convolving the images with a Gaussian with standard deviation
$\sigma(\eta)$, $G_{\sigma(\eta)} $, that is,
\begin{equation}		
I^s(\eta \mathbf{x})= G_{\sigma(\eta)} \ast I^{s-1}(\mathbf{x}),	
\end{equation}		
where we assume $\sigma(\eta)= \sigma_0 \sqrt{\eta^{-2}-1}$, $ \sigma_0 = 0.6$
and $s=1, \dots, N_{scales}$ (the number of scales).	
After the convolution, the images are sampled using bicubic                
interpolation.

We recall that the input of the algorithm is a pair of
images $I_0(\mathbf{x})$ and $I_1(\mathbf{x})$, with $\mathbf{x} = (i, j)$ the pixel
index and the output is a vector field $\mathbf{u}(\mathbf{x}) =
(u_1(\mathbf{x}), u_2(\mathbf{x}))$. The computation of
$\rho(\mathbf{u})$ involves a warping of $I_1$ and $\nabla I_1$, by the
deformation $\mathbf{u}^0$, and  the approximation field $\mathbf{u}^0$ is
computed by the multiscale scheme being zero at the coarsest
level.  A stopping criterion, for successive values of $\mathbf{u}$, $\mathbf{u^{k+1}}$ and 
$\mathbf{u^{k}}$,  is used for cessation of the algorithm before the
default number of iterations, that is,
\begin{equation}
\frac{1}{N_x N_y}\sum_{i,j}(u^{k+1}_1(i,j)-u^k_1(i,j))^2 
+ (u^{k+1}_2(i,j)-u^k_2(i,j))^2<\epsilon^2,	
\end{equation}
where $\epsilon$ is the stopping criteria threshold.

We describe the main algorithm in Algorithm 1, where the pyramidal structure is handled
and calls the function described in Algorithm 2. Algorithm 2 computes the optical flow at different scales by calling
Algorithm 3 or Algorithm 4. Algorithm 3 is
the function {\it split-Bregman-1}, related to the model that has the gradient as the regularising term and  Algorithm 4 
is the function  {\it split-Bregman-alpha}, related to the model 
that considers the fractional operator as the regularising term.
The method that calls the algorithm {\it split-Bregman-1}  is hereafter denoted by TV-L1-SB method and the
one that calls the algorithm {\it split-Bregman-alpha} named the TV-L1-SB-alpha method.

\begin{algorithm}[h]

\caption{{\it Main Algorithm} -- {\it Create the Pyramidal Structure}}
				 
\hspace{0.5cm}{\bf Input:}  $I_0, I_1$, $\lambda, \theta, \epsilon, \eta, N_{scales},
N_{warps}, N_{maxitr}$
				  
\hspace{0.5cm}{\bf Output:}  Optical flow $\mathbf{u} =(u_1, u_2)$.

\hspace{0.5cm}Normalize images between $0$ and $255$ and convolve the images with a Gaussian of $\sigma_0 =0.6$
				 

\hspace{0.5cm}Create the pyramid of images $I^s$ using $\eta$ (with $s= 1, ..., N_{scales}$)

\hspace{0.5cm}$(u_1^{N_{scales}},u_2^{N_{scales}})=(0,0)$
				 
\hspace{0.5cm}{\bf for } $s = N_{scales}$ {\bf to} $1$ {\bf do}
				 
\hspace{1cm} {\it Optical-Flow-SB}
($I_0,I_1, \mathbf{u}^0, \theta, \lambda_{\ell}, \epsilon, N_{warps}, N_{maxitr} $)
				 
\hspace{1cm} {\bf if} $s>1$ {\bf then}
		 
\hspace{1.5cm}$u^{s-1}(\mathbf{x})= \frac{1}{\eta} \mathbf{u}^s(\mathbf{x}/\eta)$
				 
\hspace{1cm}{\bf end}
				 
{\hspace{0.5cm}\bf end}
\end{algorithm}
\begin{algorithm}[h]

\caption{{\it Optical-Flow-SB} ($I_0,I_1, \mathbf{u}^0, \theta,
\lambda_{\ell}, \epsilon, N_{warps}, N_{maxitr} $)}

Compute the optical flow at different scales
				  				  
\hspace{0.5cm}{\bf for} $w = 1$ {\bf to} $N_{warps}$ {\bf do}
				  
\hspace{0.5cm}Compute $I_1(\mathbf{x}+\mathbf{u}^0), \nabla I_1(\mathbf{x}+\mathbf{u}^0)$ using bicubic interpolation
				  
				   
\hspace{1cm}{\bf while} $n<N_{maxitr}$ {\bf and} stopping\_criterion $>\epsilon$ {\bf do}
				 
\hspace{1.5cm}$\mathbf{v}= \mathbf{u} + TH(\mathbf{u},\mathbf{u}^0)$
				  
\hspace{1.5cm}$u_1=$ {\it split-Bregman-alpha} $(v_1, \theta, \lambda_{1},\alpha)$ or
{\it split-Bregman-1} $(v_1, \theta, \lambda_{1})$
				  
\hspace{1.5cm}$u_2=$ {\it split-Bregman-alpha} $(v_2, \theta, \lambda_{2},\alpha)$ or
{\it split-Bregman-1} $(v_2, \theta, \lambda_{2})$
				 
				  
\hspace{1cm}{\bf end}
				 
\hspace{0.5cm}{\bf end}
\end{algorithm}
\begin{algorithm}[h]

\caption{{\it split-Bregman-1}  $( v, \theta, \lambda_{\ell})$}

Compute the split-Bregman iteration for the gradient regularization term
				
\hspace{0.5cm}Fix  TOL  and  $u^0=v,$ $d^0= b^0=0$
				
\hspace{0.5cm}{\bf while} $|u^k-u^{k-1} |>$TOL {\bf do}
				
\hspace{1cm} $u^{k+1} =  \displaystyle  \frac{1}{\frac{1}{\theta} 
+ 4\lambda_{\ell}}\Big[\lambda_{\ell} U^k  + \frac{1}{\theta} v 
-   \lambda_{\ell}  \mbox{div}( d^k  - b^k)\Big]$
				
\hspace{1cm}$ \frac{\partial  u^{k+1}}{\partial n}  = (d^k  - b^k) \cdot n \quad $  in $\partial \Omega$,
				
\hspace{1cm}$d^{k+1} = shrink\left(\nabla u^{k+1} + b^k, \frac{1}{\lambda_{SB}}\right)$,
				
\hspace{1cm}$b^{k+1}=b^k + \nabla u^{k+1} - d^{k+1}$.
				
\hspace{0.5cm}{\bf end}
\end{algorithm}

\begin{algorithm}[h]

\caption{{\it split-Bregman-alpha} $ (v, \theta, \lambda_{\ell}, \alpha)$}

Compute the split-Bregman iteration for the fractional regularization term
		
\hspace{0.5cm}Fix  TOL  and  $u^0=v,$ $d^0= b^0=0$
				
\hspace{0.5cm}{\bf while} $|u^k-u^{k-1} |>$ TOL {\bf do}
				
\hspace{1cm}
$\ds{u_{i,j}^{k+1}=
\frac{\lambda_{\ell}}{a_\ell}
\left[\sum_{p=i}^{N_x}w^{(\alpha)}_{p-i}\sum_{\ell=0,\ell\neq i}^p w^{(\alpha)}_\ell u^k_{\ell-i,j}
+\sum_{p=j}^{N_y} w^{(\alpha)}_{p-j}\sum_{\ell=0,\ell\neq j}^p w^{(\alpha)}_\ell
u^k_{i,\ell-j}\right.}
$

\hspace{1.5cm}$\ds{\left. +\frac{1}{\theta\lambda_{\ell}} (v_\ell)_{i,j}+\sum_{p=0}^{N_x-i}w^{(\alpha)}_{p}(d_{1}^k - b_{1}^k)_{i+p,j} 
+ \sum_{p=0}^{N_y-j}w^{(\alpha)}_{p}(d_{2}^k  - b_{2}^k)_{i,j+p}\right]},		
$	
				
\hspace{1cm}$ u^{k+1}  = 0  \quad $  in $\partial \Omega$
				
\hspace{1cm}$d^{k+1} = shrink\left(D_{-}^{\alpha} u^{k+1} + b^k, \frac{1}{\lambda_{SB}}\right)$
				
\hspace{1cm}$b^{k+1}=b^k + D_{-}^{\alpha} u^{k+1} - d^{k+1}$
				
\hspace{0.5cm}{\bf end}
\end{algorithm}

The algorithm depends on several parameters. Therefore for a better understanding of
the algorithm written in this section, we give a brief overview of the
role of each parameter.

The {\it data attachment weight} $\lambda$ is a relevant parameter
and it determines the smoothness of the output. The smaller this parameter is,
the smoother the solutions we obtain.  It depends on the range of motions of
the images, so its value should be adapted to each image sequence.
The {\it tightness} $\theta$ serves as a link between the attachment and the
regularization terms. It should have a small value in order to
maintain both parts in correspondence. 
The {\it penalty parameters in split Bregman iteration} 
$\lambda_{\ell}, \ \ell=1,2$ are used in order to holding the penalty constraint.
The choice of these parameters are related with the value choosen for
$\theta$ in order to reach a faster convergence.
The {\it stopping criteria threshold} $\epsilon$ 
is a trade-off between precision and running time. 
A small value will yield more accurate solutions at the expense 
of a slower convergence.
The {\it downsampling factor} $\eta$ is used  in order to downscale the
original images to create the pyramidal structure and ranges values between
$0$ and $1$.
The {\it number of scales} $N_{scales}$
is used to create the pyramid of images.  
The {\it number of warps} $N_{warps}$			 
represents the number of times that $I_1(\mathbf{x}+\mathbf{u}_0)$ 
and $\nabla I_1(\mathbf{x}+\mathbf{u}_0)$ are computed per scale.			 
It  affects the running time.
The {\it parameter}  $\alpha$ {\it of the fractional operator}
affects the regularisation operator and ranges values between $0$ and $2$.

\section{Experimental results}

The numerical method has been implemented and applied to different
image sequences and all the results presented here uses  two frames of a sequence
of images as input.
First we discuss the TV-L1-SB method by comparing the resulting computing flow
with the one obtained by a method given in \cite{san2013} and,
based in an evaluation methodology well established by now, we notice
that in general the TV-L1-SB method performs better.
The second section discusses
the effect of the fractional regularisation operator. 
It is shown  that the parameter $\alpha$ should be adjusted depending
on the geometry or texture complexity of the various regions of the
 image. 
 
There are several parameters involved in the algorithm, that we list 
in Table 1 with a short explanation and the
range of the best values  discussed in the next section.
\begin{table}[h]
\centerline{
\begin{tabular}{|l  l l|}\hline
{ Parameter} & Description & Range of values\\\hline
$\lambda $ & data attachment weight& [0.1 \  1]\\
 $\theta$ & tightness & [0.1 \  1]\\
 $\epsilon$ & stopping threshold & 0.01 \\
 $\eta$ & zoom factor & 0.5\\
 $N_{scales}$ & number of scales & [3 \ 5]\\
 $N_{warps}$ & number of warps & [3 \ 5]\\
 $\lambda_{SB}$ & split-Bregman parameter $\lambda_{SB}=\lambda_1=\lambda_2$ &
 [1 \ 10]\\
 $\alpha$ & order of fractional operator & [0 \ 2]\\\hline
\end{tabular}}
\caption{Parameters involved in  the numerical methods}
\end{table}

\subsection{Performance of the TV-L1-SB method}

In this section the goal is  to show the advantage of applying the split Bregman technique
in the determination of the optical flow.  To that end, we consider a very recent
numerical method, used to solve the same
optical flow model, presented in  \cite{san2013}. The method of \cite{san2013} is hereafter 
 named the TV-L1 method and we compare its performance with the TV-L1-SB method.
 The  TV-L1-SB method, described in
section 3.1 and 3.3, differs essentially
from the TV-L1 method in the application of
the split Bregman method  to solve the first minimisation problem (\ref{first_minimization}).

To show the performance of both numerical methods we use a test sequence of
different images  from the Middlebury database, displayed in Figure 1.


\begin{figure}[h]
\centerline{
\includegraphics[width=4cm,height=3.2cm]{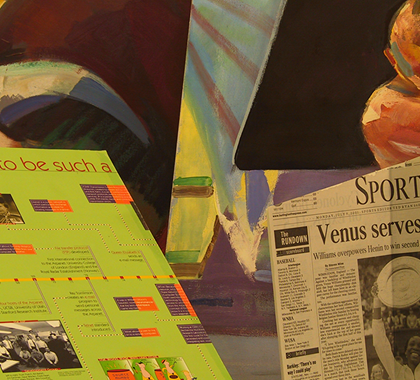}	
\includegraphics[width=4cm,height=3.2cm]{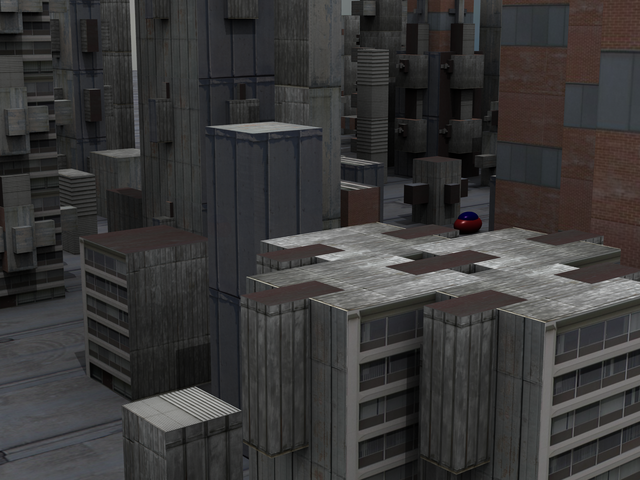} 
\includegraphics[width=4cm,height=3.2cm]{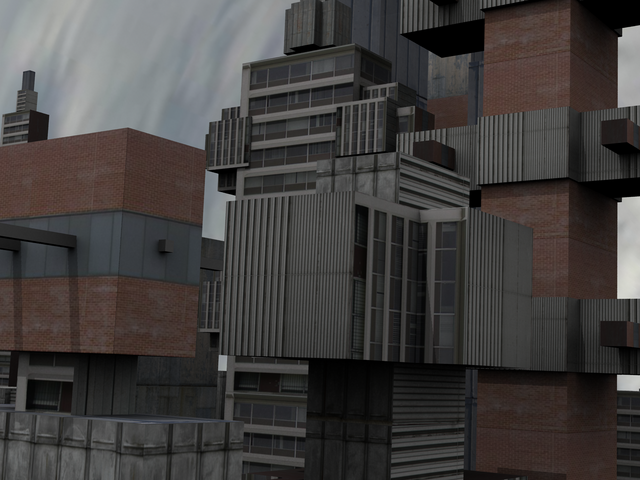} 
\includegraphics[width=4cm,height=3.2cm]{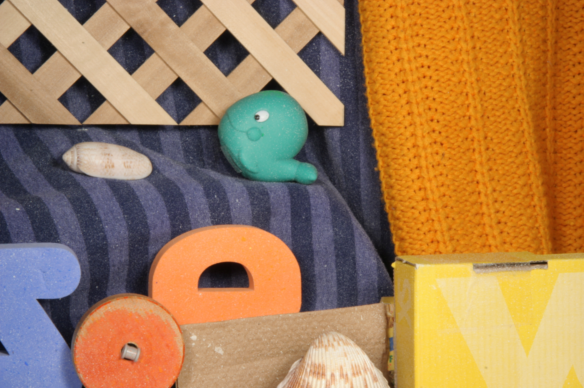}}
\vspace*{0.1cm}
\centerline{
\includegraphics[width=4cm,height=3.2cm]{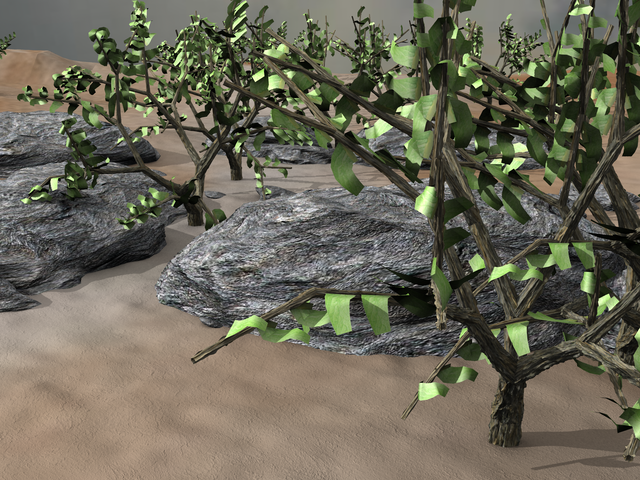} 
\includegraphics[width=4cm,height=3.2cm]{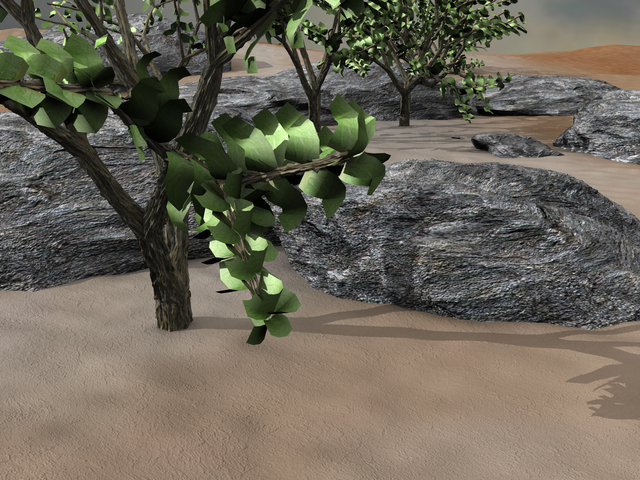} 
\includegraphics[width=4cm,height=3.2cm]{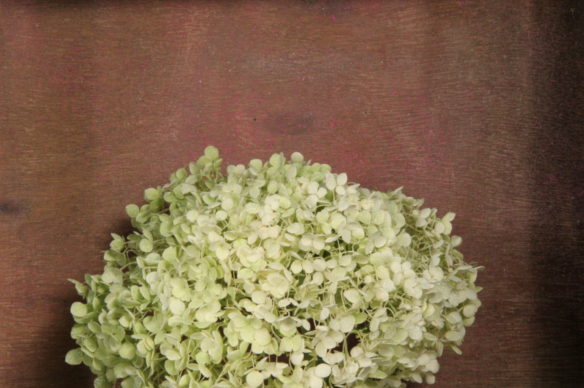} 
\includegraphics[width=4cm,height=3.2cm]{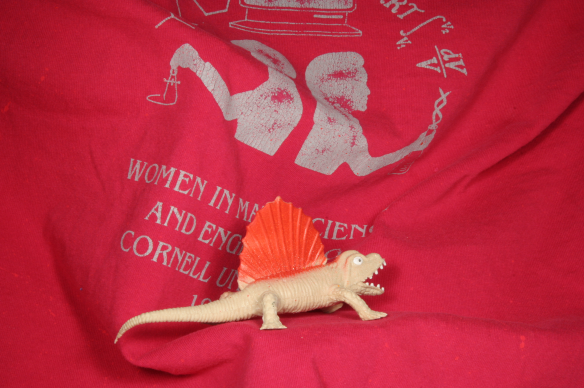}}
\caption{Image sequences from the Middlebury database. One frame of each sequence is displayed.
Top and from left to right: Venus, Urban2, Urban3,
Rubber Whale; Bottom and  from left to right: Grove3, Grove2, Hydrangea, Dimetrodon.}
\end{figure}

The most common
used measure of performance for optical flow is the angular error (AE) between a
flow vector $(u_1,u_2)$ and the ground truth flow $(gt_1,gt_2)$.  It is the
angle in $3D$ between $(u_1,u_2,1)$ and $(gt_1, gt_2, 1)$.
The AE is calculated by the following formula:		
\begin{equation}
{\rm AE} = {\arccos}\left( \frac{1+ u_1\times gt_1 + u_2 \times gt_2}
{\sqrt{1+ u_1^2 + u_2^2}\sqrt{1+gt_1^2+gt_2^2}}\right).
\label{AE}	
\end{equation}		
The popularity of this measure is due to the seminal work  by Barron et al \cite{bar1994}.
This measure provides a relative measure of performance that avoids the divide
by zero problem for zero flows. Errors in large flows are penalized less in AE
than errors in small flows. 
The AE also contains an arbitrary scaling constant $1$ to convert the units 
from pixels to degrees.

A complementar measure is usually used, which is,
the End Point Error (EPE) defined by	
\begin{equation}		
{\rm EPE} =\sqrt{(u_1-gt_1)^2 + (u_2-gt_2)^2}.		
\label{EPE}
\end{equation}
This measure may be more appropriate for some applications.
Therefore we report both herein.

These measures are computed at each pixel and consequently the corresponding averages
of AE (AAE) and of EPE (AEPE) are used.
In some cases the standard deviations of AE (SDAE) will be also presented.

In the next experiments some of the values of the parameters in our model are
chosen to be the values for which the TV-L1 method performs better, according to \cite{san2013},
in order to compare its performance with the performance of the  TV-L1-SB method.
Regarding the choice of the parameters $\lambda_\ell$, $\ell=1,2$,
associated with split Bregman method,
we assume 
$$\lambda_{SB}=\lambda_1=\lambda_2.$$
Note  that the parameters $\lambda_\ell$, $\ell=1,2$ are related, respectively, with the
flow velocity $u_\ell, \ell=1,2$.
In \cite{gol2009} the authors have found  that for a faster convergence
a good choice for the split Bregman parameter, $\lambda_{SB}$,
can be
$
\lambda_{SB}= {2}/{\theta}.
$
For the TV-L1-SB method this seems to be also a  suitable choice.

We start with some discussion about the choice of the split Bregman parameter $\lambda_{SB}$.
Figure 2 displays the results of the experiments done 
for the Rubber-Whale sequence (see Figure 1). It shows the performance of the TV-L1-SB method for
different values of  $\lambda_{SB}$ when $\theta=0.4$
and for different data attachment weights $\lambda$.
The tightness parameter $\theta=0.4$ is chosen according to \cite{san2013} as mentioned previously.
It can be seen that the TV-L1-SB method performs better for the set of split Bregman values, $1\leq \lambda_{SB}\leq 10$,
reaching its best between $3$ and $10$.

\begin{figure}[h]
\centerline{\includegraphics[height=50mm,width=75mm]
{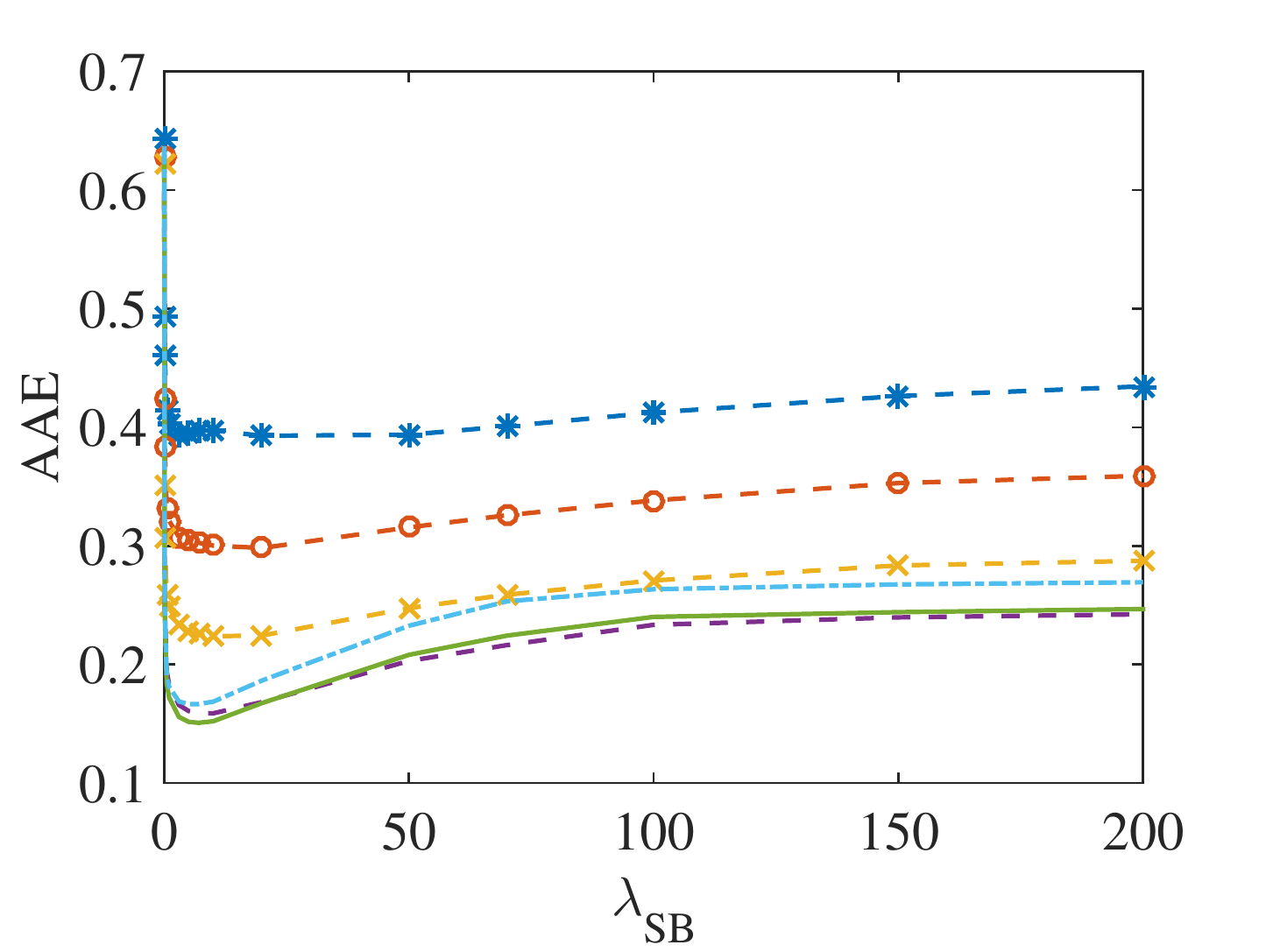}
\includegraphics[height=50mm,width=75mm]
{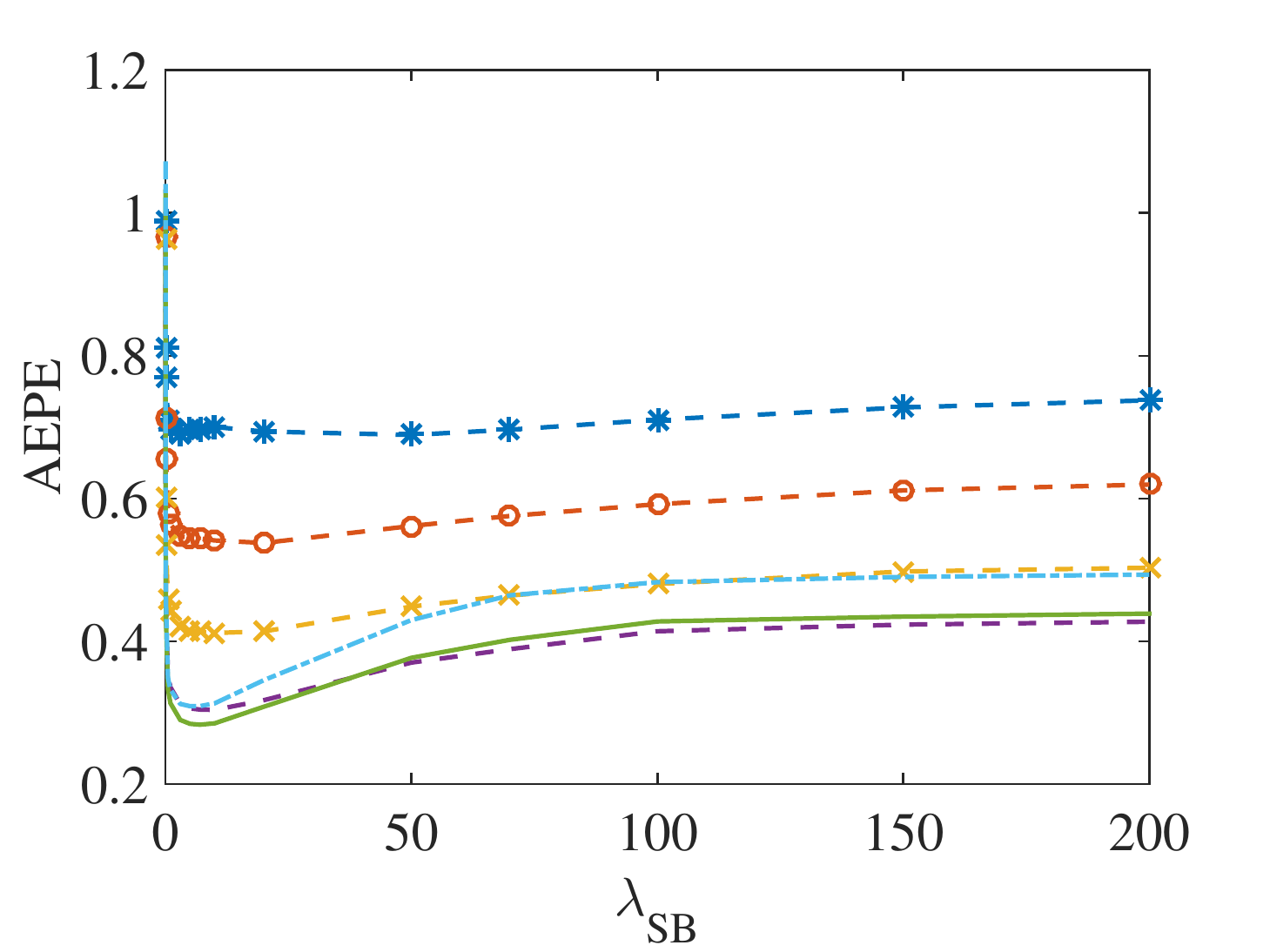}}
\caption{Performance of the TV-L1-SB method for the Rubber Whale sequence:  $\theta=0.4$ and
   $Nwraps=Nscales=5$. The method performs better for values of the split Bregman parameter
   $3\leq\lambda_{SB}\leq 10$. Plots for different values of the data attachment weight $\lambda$: 
 $\lambda=0.03 (\-- *)$; $\lambda=0.05(\-- o)$; $\lambda=0.1(\-- \times)$; 
 $\lambda=0.3(\-- \--)$; $\lambda=0.5(\--)$;
$\lambda=1(\-- \cdot \--)$.
}
\label{fig2}   
\end{figure}

In Table \ref{tvl1table} we  compare the TV-L1-SB method
with the  TV-L1 method  for the image sequences presented in Figure 1. 
For each image sequence we have chosen the parameters $\lambda$ (data attachment weight) and $\theta$
(tightness) for which the TV-L1 method performs better. The split Bregman parameter $\lambda_{SB}$ 
needed in the TV-L1-SB method is assumed to be $10$ for all image sequences. The results point out
the TV-L1-SB method performs better than the TV-L1 method. 
\begin{table}[h]
\centerline{	
\begin{tabular}{|llll|l||llll|l|}\hline
Method & AAE & AEPE & SDAE & 	Data & Method & AAE & AEPE & SDAE & Data \\\hline		
Grove2 &          &            &              & 6 scales & \multicolumn{2}{ l }{RubberWhale} &&& 4 scales\\	
TV-L1  &  0.7070 &  2.1876  & 0.4351  & $\lambda= 0.3$ &	TV-L1&  0.2281 &  0.4155 &  0.2724 & $\lambda= 0.4$\\		
TV-L1-SB & 0.5443  & 1.7974  &  0.4283 & $\theta = 0.3$ & TV-L1-SB &  0.1530 & 0.2905  &   0.2406 & $\theta = 0.4$ \\\hline
Grove3&&&& 4 scales &  \multicolumn{2}{ l }{Hydrangea}&&& 4 scales\\
TV-L1  &  0.5856 & 2.8451 & 0.4803 & $\lambda= 0.5$ & TV-L1&  0.4185 & 2.1619 & 0.2798 & $\lambda= 0.1$ \\	
TV-L1-SB & 0.4198  & 2.4920 & 0.3984  & $\theta = 0.4$ & TV-L1-SB &   0.3107 & 1.9487  & 0.2049  & $\theta = 0.8$\\\hline
Urban2 &&&& 6 scales & \multicolumn{2}{ l }{Dimetrodon}&&& 5 scales\\	
TV-L1 &  0.6623 &  7.4200 &  0.5477 & $\lambda= 0.5$ & TV-L1 &  0.4282 & 1.1511 &   0.3109 & $\lambda= 0.3$\\		
TV-L1-SB&  0.5319  & 7.1988  &  0.4941  & $\theta = 0.3$ & TV-L1-SB&   0.3356 &  1.0051  & 0.2724  & $\theta = 0.3$ \\\hline	
Urban3&&&& 5 scales & Venus&&&& 4 scales\\
TV-L1 &  0.9403 &  6.4110  &  0.6624 & $\lambda= 0.9$ & TV-L1 &  0.6693 &  2.8136  & 0.4318  & $\lambda= 0.4$\\
TV-L1-SB & 0.7890  & 5.9711  & 0.6606  & $\theta = 0.7$ & TV-L1-SB &  0.4629 &  2.4105 &  0.3757 & $\theta = 0.6$\\\hline	
\end{tabular}}		
\caption{ Errors AAE, AEPE and SDAE for the TV-L1 method
and the TV-L1-SB method for the sequence of images presented in Figure 1.}
\label{tvl1table}		
\end{table}

\begin{figure}[h]
\centerline{\includegraphics[height=50mm,width=75mm]
{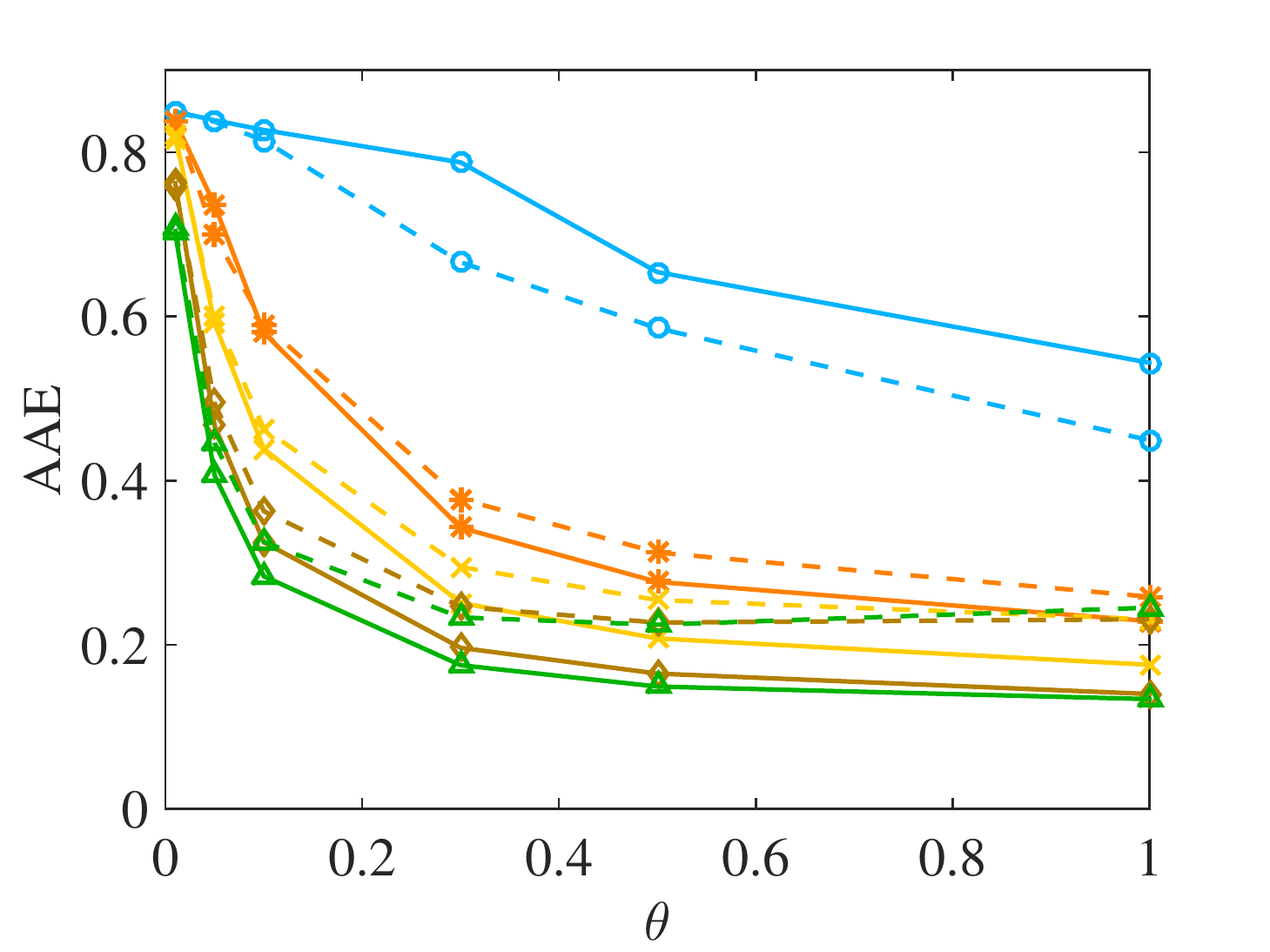}
\includegraphics[height=50mm,width=75mm]
{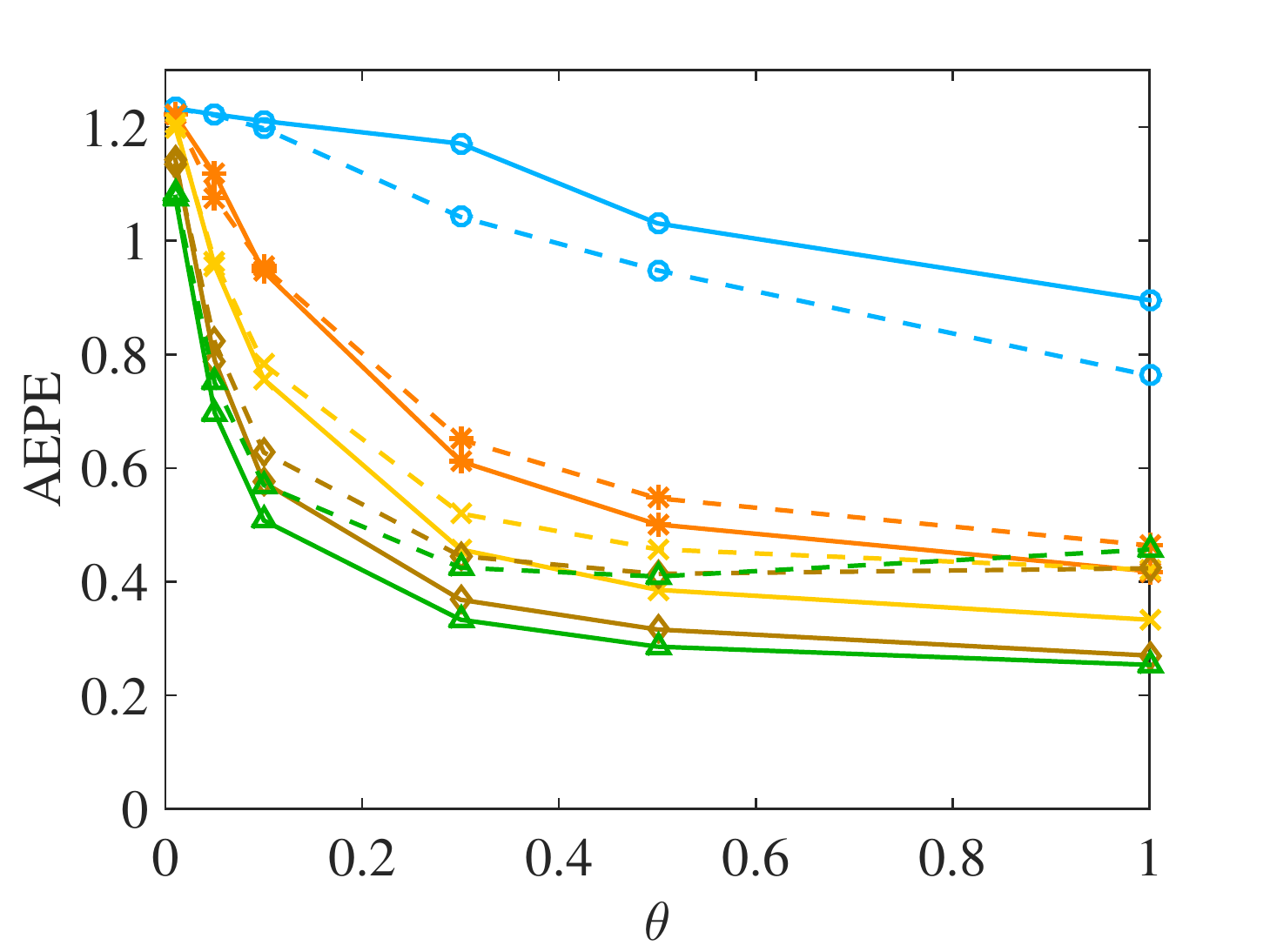}}
\caption{Errors AAE and AEPE for the Rubber Whale sequence.  The tightness parameter $\theta$ changes and
$Nwraps=5$, $Nscales=4$. The TV-L1 method is represented by the dashed  lines $(\-- \--)$
and the TV-L1-SB method by the solid lines $(\--)$.
The different values of  the data attachment weight $\lambda$ are marked by different symbols: 
$\lambda=0.01 ( o )$, $\lambda=0.05 ( * )$, $\lambda=0.1 ( \times )$, $\lambda=0.2
( \diamond )$,
$\lambda=0.3 ( \ \hat{ } \ )$.}
\label{figtheta}   
\end{figure}
\begin{figure}[h]
\centerline{\includegraphics[height=50mm,width=75mm]
{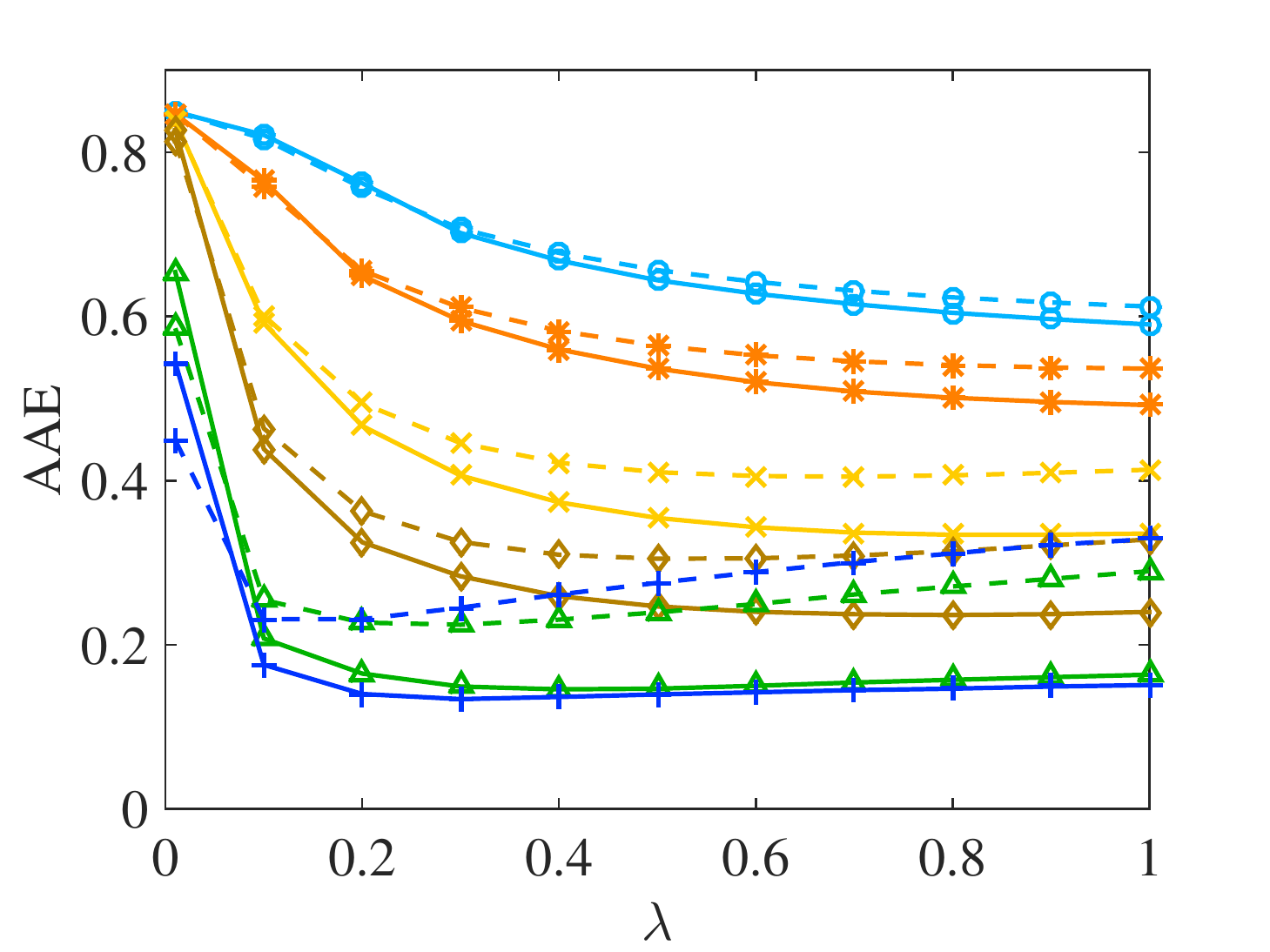}
\includegraphics[height=50mm,width=75mm]
{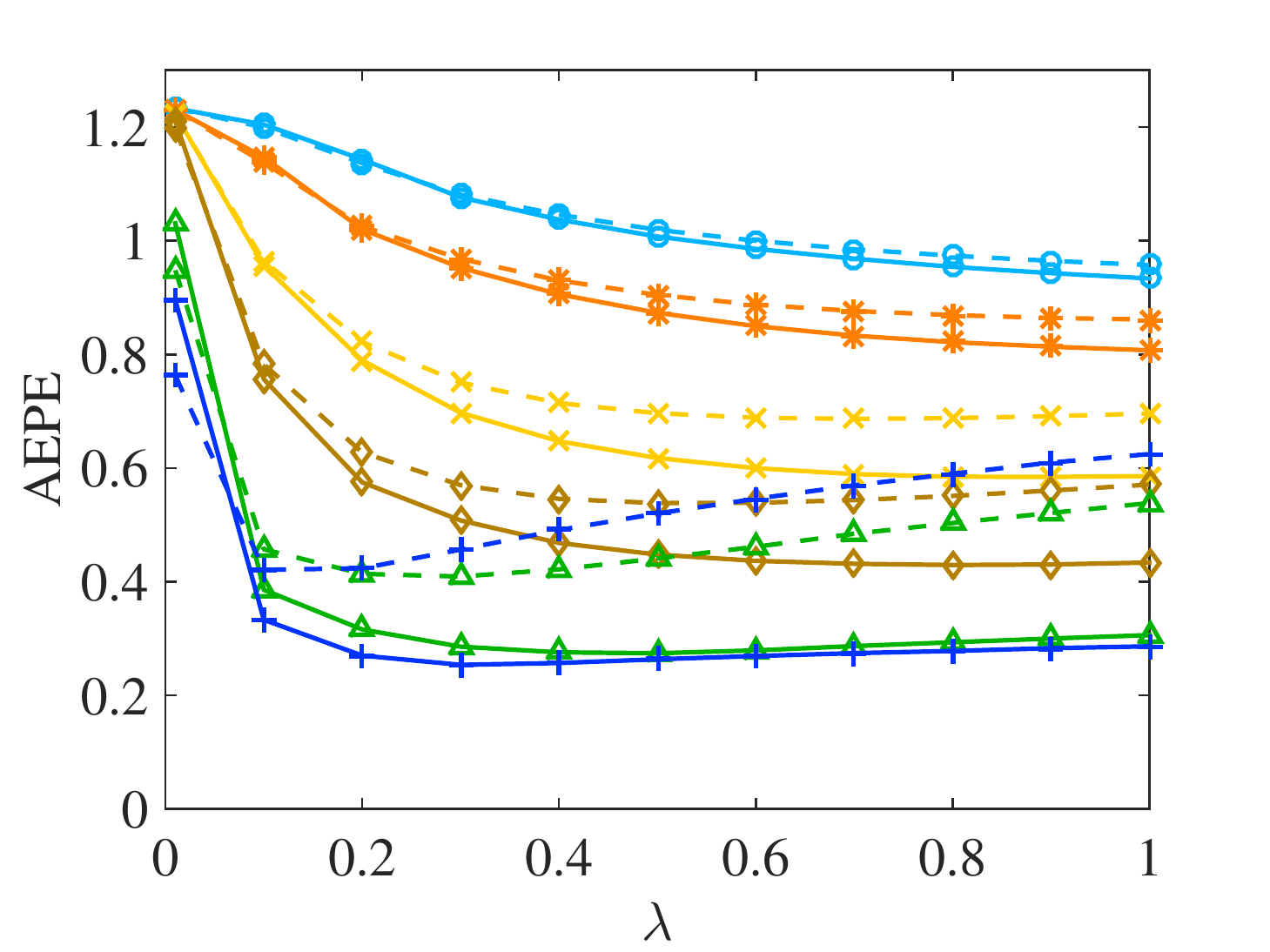}}
\caption{Errors AAE and AEPE for the Rubber Whale sequence: The data attachment weight parameter
$\lambda$ changes and
  $Nwraps=5$, $Nscales=4$. The TV-L1 method is represented by  the dashed lines $(\-- \--)$
and the TV-L1-SB method by the solid lines $(\--)$.
The different values of the tightness parameter $\theta$ are marked with different symbols: 
$\theta=0.01 ( o )$, $\theta=0.02 ( * )$, $\theta=0.05 ( \times )$, $\theta=0.1
( \diamond )$,
$\theta=0.5 ( \ \hat{ } \ )$, $\theta=1 ( \ + \ )$.}
\label{figlambda}   
\end{figure}

To have a more complete perspective between the differences on the performance
of both methods we exhibit in Figures \ref{figtheta}  and  \ref{figlambda}  additional results,
with different values of $\theta$
and $\lambda$ in the Rubber Whale case and for $\lambda_{SB}=10$.
The results confirm the TV-L1-SB method presents 
smaller errors in general.
By running these experiments we have observed that the TV-L1-SB method
is slower if the split Bregman parameter $\lambda_{SB}$ is very far away from the estimate $2/\theta$.
We have executed the algorithm for a fixed $\lambda_{SB}=10$ and as $\theta$ becomes
larger the method becomes slower, suggesting the parameter 
$\lambda_{SB}$ should  be adjusted depending on $\theta$, for a faster
convergence. We did not adjust the parameter $\lambda_{SB}$ to plot  Figures
\ref{figtheta}  and  \ref{figlambda}  to emphasize the effect of only changing the
 data attachment weight  parameter $\lambda$ and tightness parameter $\theta$.

\subsection{The effect of the fractional regularization operator}

The main purpose of this section is to
present the effect of the parameter $\alpha$ in the
estimation of the optical flow. Only partial regions of the image are used in
order to show efficiently the influence of $\alpha$.
Three main regions are considered:  with edges, corners and flat
that can also present high or low texture and motion discontinuities.


\begin{figure}[h]			
\includegraphics[width=4cm,height=3.2cm]
{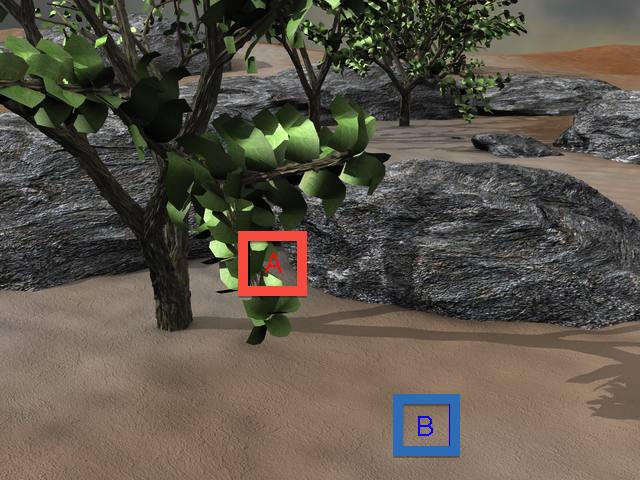} 
\includegraphics[width=4cm,height=3.2cm]
{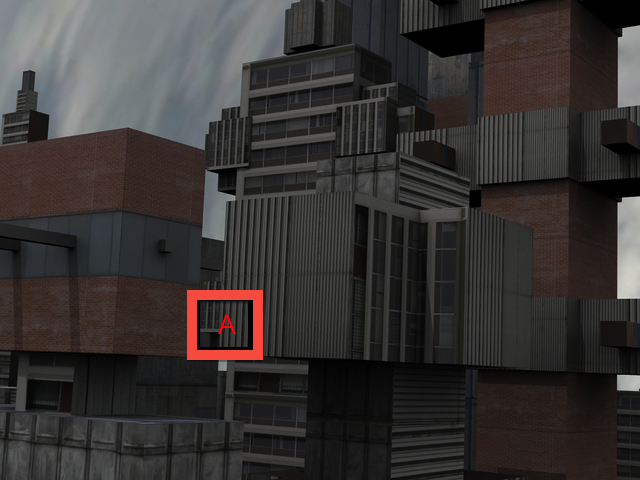} 
\includegraphics[width=4cm,height=3.2cm]
{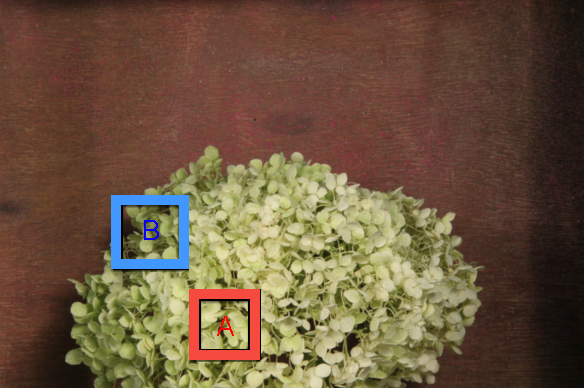} 
\includegraphics[width=4cm,height=3.2cm] 	
{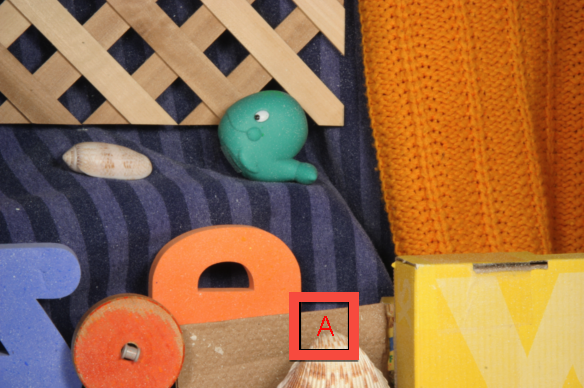}
\caption{The regions considered for each image sequence are marked as regions A (red)
or B (blue). From  left to right: Grove2 (Regions A and B), Urban3 (Region A),
Hydrangea (Regions A and B) and Rubber Whale (Region A).}			
\label{fracimage}		
\end{figure}

We start to discuss the results for the image sequences considered the most difficult
from the Middleburry database, which are Grove and Urban. For these data images
the optical flow  methods usually perform worst. As described in \cite{bak}
Grove contains a close up view of a tree, with a substantial parallax and motion
discontinuities and Urban contains the image of a city with substantial 
motion discontinuities, a large motion range and an  independently moving object.
We consider the image sequences Grove2 and Urban3 displayed in
Figure \ref{fracimage}. 

\begin{figure}[h]
\includegraphics[height=40mm,width=56mm]
{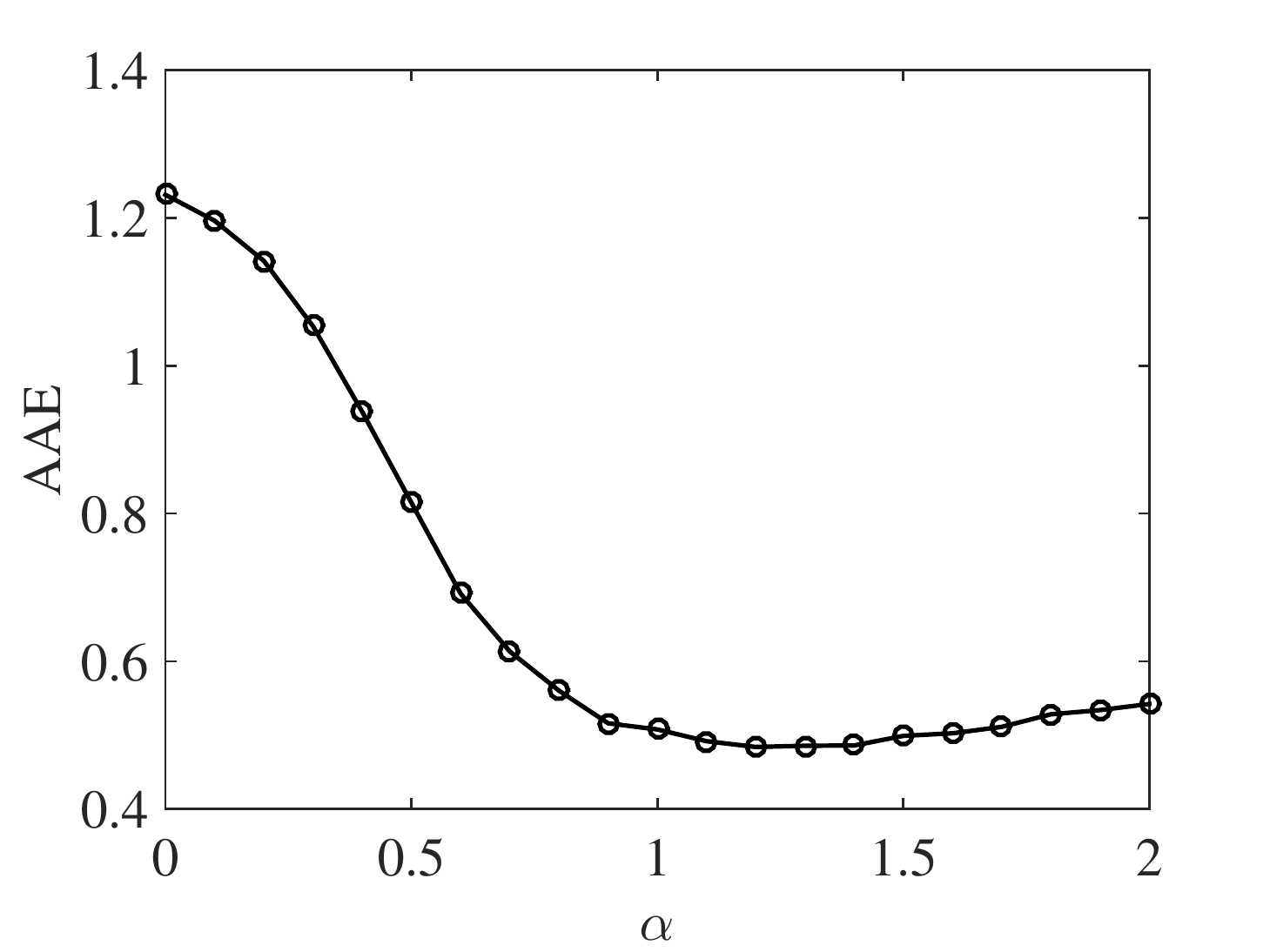}
\includegraphics[height=40mm,width=56mm]
{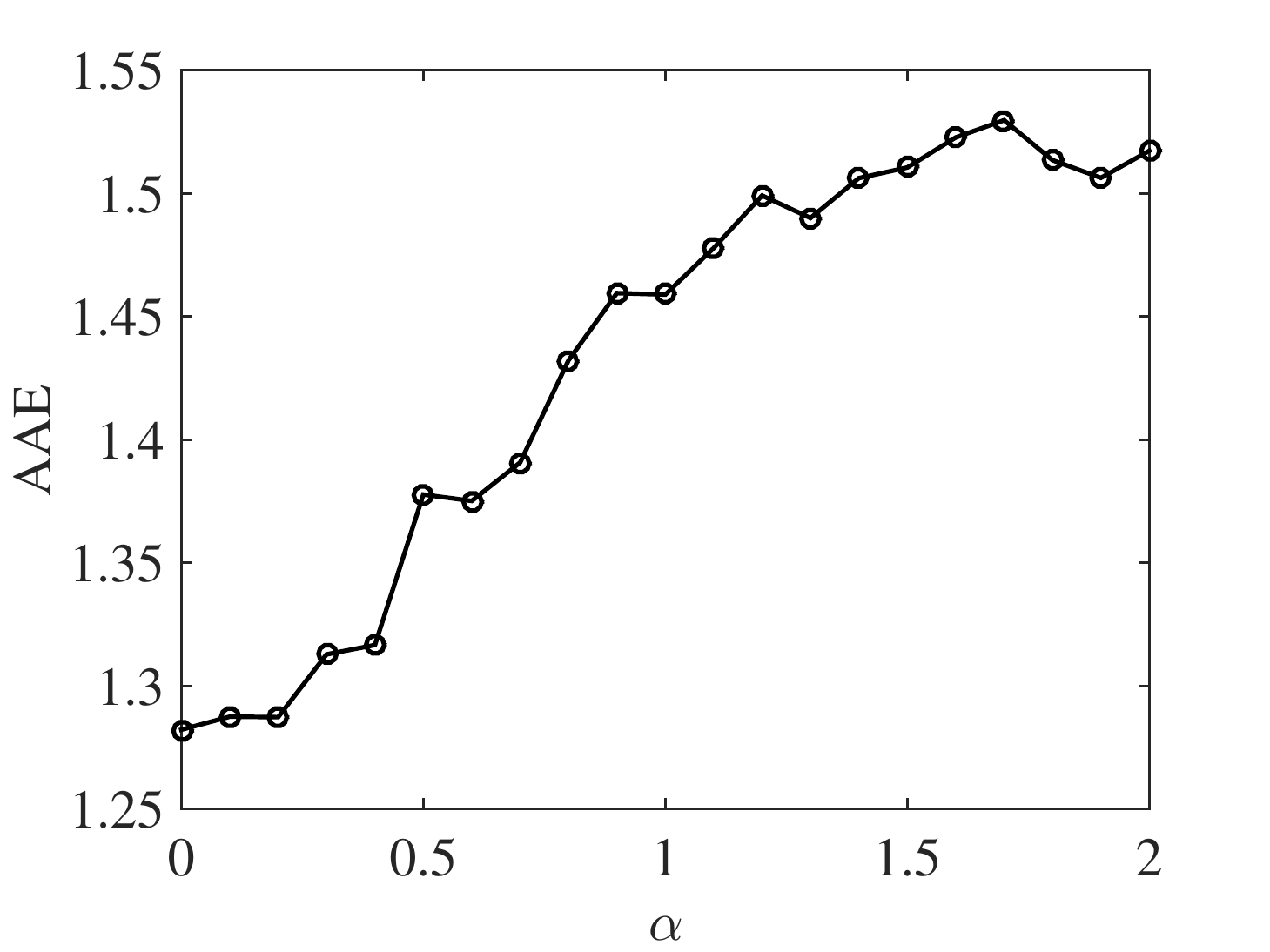}
\includegraphics[height=40mm,width=56mm]
{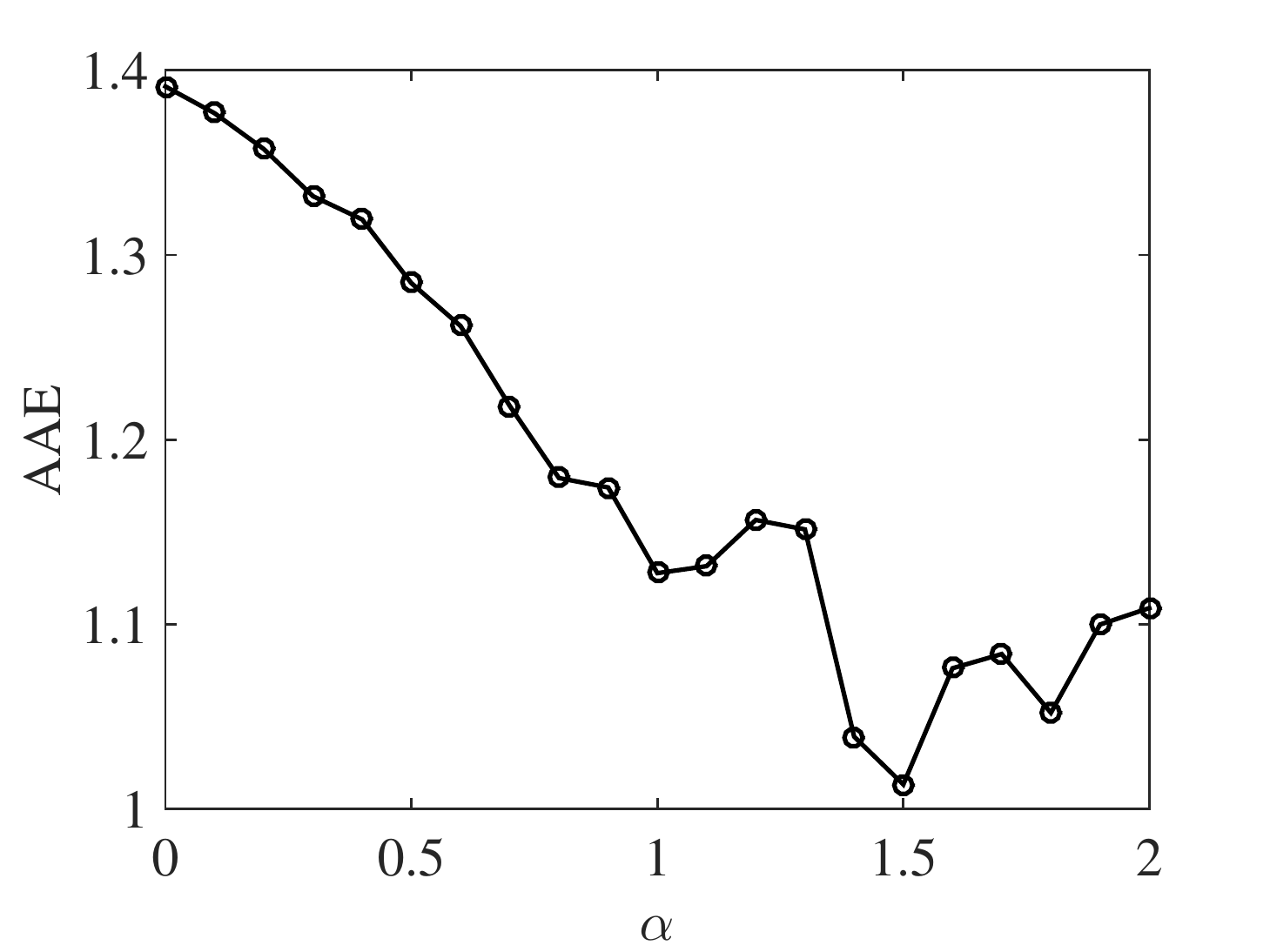}
\newline
\includegraphics[height=40mm,width=56mm]
{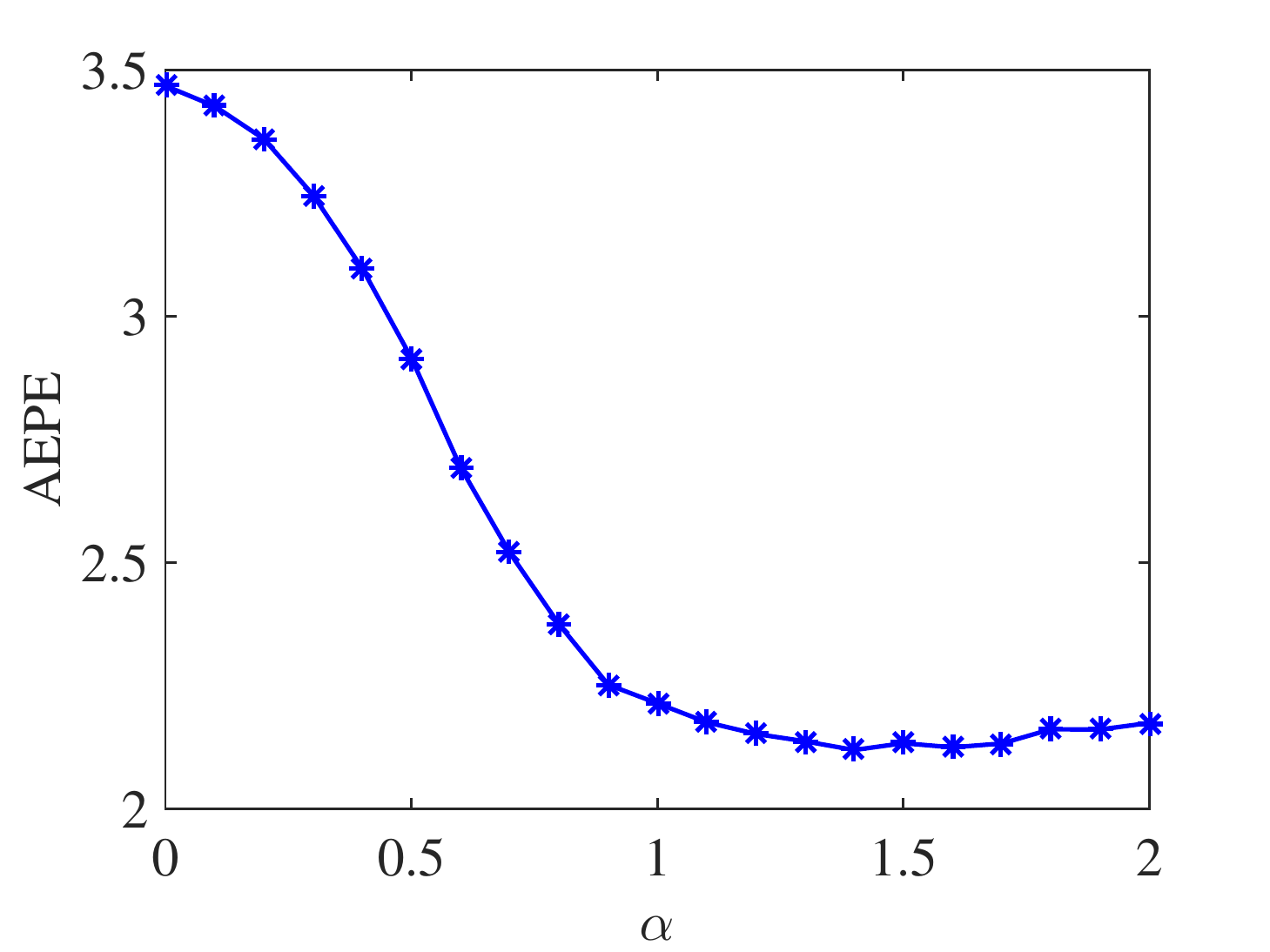}
\includegraphics[height=40mm,width=56mm]
{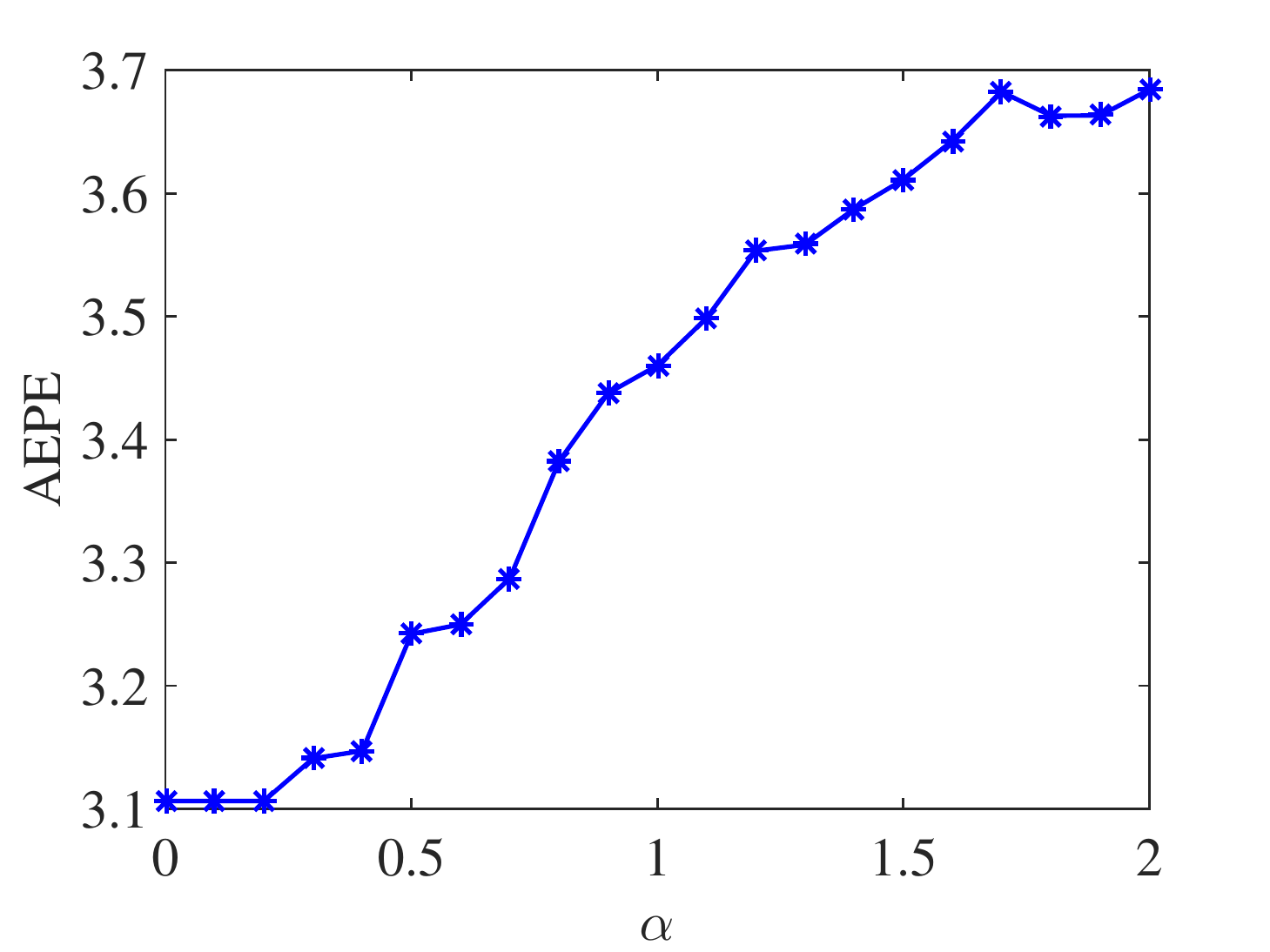}
\includegraphics[height=40mm,width=56mm]
{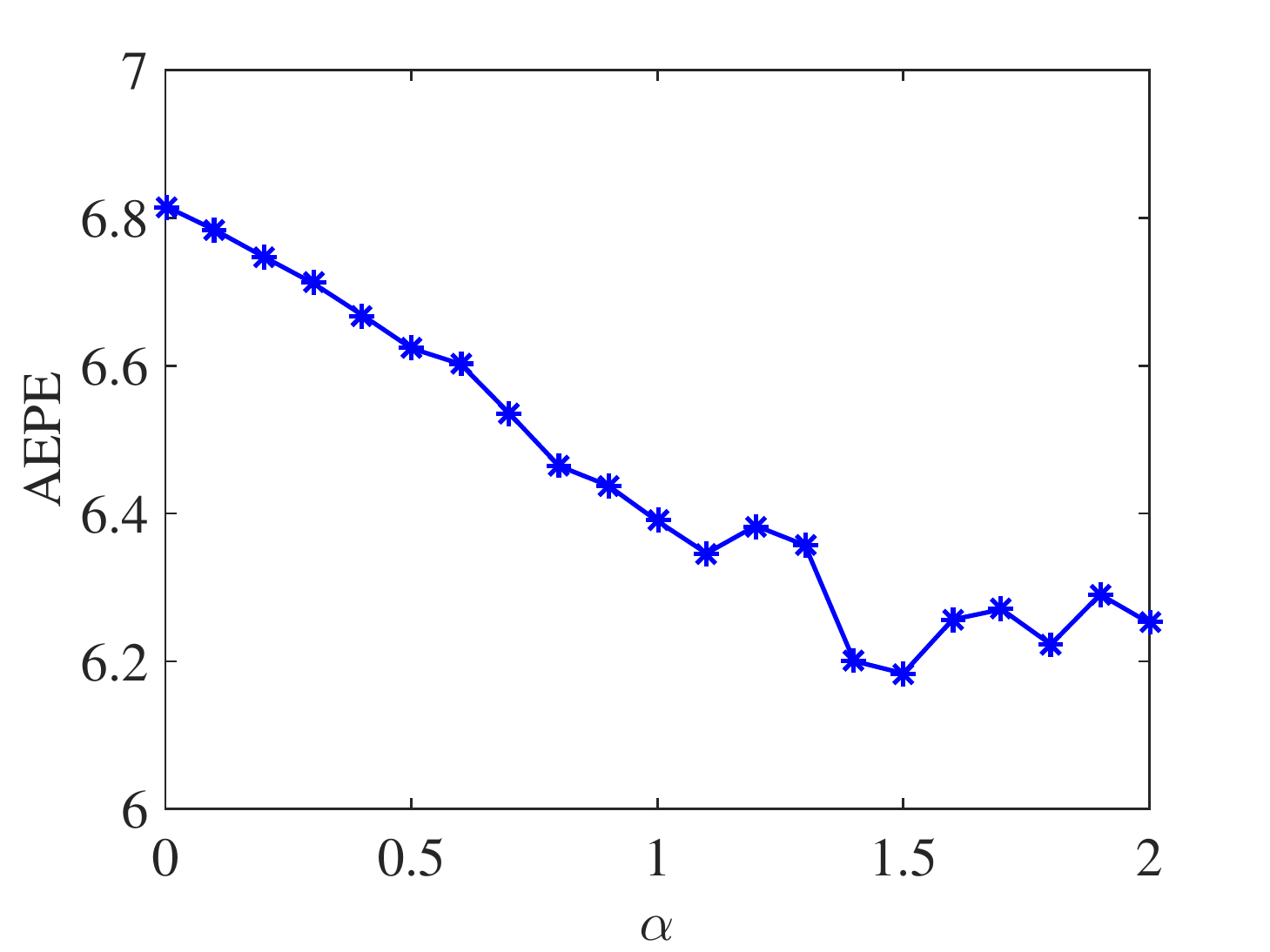}
{\small \hspace*{2.6cm} (a) \hspace*{5.2cm}(b)  \hspace*{5.2cm}(c)}
\caption{(a) Grove2 -- Part A. For  $\lambda= 0.3$ and $\theta = 0.3$.
Best results for $\lambda_{SB}=5$. AAE: best result
  for $\alpha=1.2$. AEPE: best result for $\alpha=1.4$;
 (b) Grove2 -- Part B. For  $\lambda= 0.3$ and $\theta = 0.3$.
Best results for $\lambda_{SB}=7$. AAE: best result
  for $\alpha=0$. AEPE: best result for $\alpha=0$;
 (c) Urban3 -- Part A. For  $\lambda= 0.9$ and $\theta = 0.7$.
Best results for $\lambda_{SB}=1$. AAE: best result
  for $\alpha=1.5$. AEPE: best result for $\alpha=1.5$.}
\label{Grove2Urban3}   
\end{figure}

For the
sequence Grove2 we have selected two types of regions for discussion: one represents 
a flat region (marked as region B)
and the other one represents a region with edges and a parallax effect (marked as region A).
We assume $N_{wraps}=5$, $N_{scales}=5$ and the tightness parameter
$\theta$ and the data attachment weight parameter $\lambda$ are chosen according to Table 2, that is,
$\theta=0.3$ and $\lambda=0.3$.
In Figures \ref{Grove2Urban3}(a)  and \ref{Grove2Urban3}(b) we report the errors for the optical flow computed in the two regions
 of Grove2 for different values of 
$\alpha$. We have also done tests for different values of the  split Bregman 
parameter,  $\lambda_{SB}$, in particular with values between $1$ and $10$. In Figure \ref{Grove2Urban3}(a),
related to the region A marked in Figure \ref{fracimage},
we only plot the results for the split Bregman parameter $\lambda_{SB}=5$,  the value for which
we have obtained the best results. It is shown
the best result, represented by the smallest errors, have been reached  
at $\alpha=1.2$ for the error AAE and at $\alpha=1.4$ for the error
AEPE. In general, the best  results are for values of $\alpha$ between $1$ and $1.5$. 
In Figure \ref{Grove2Urban3}(b), we show what happens for the flat region, 
region B of Grove2.  The best results in this case have been obtained for
$\lambda_{SB}=7$ and the errors AAE and AEPE
are smaller  for $\alpha=0$.

We turn now to the image sequence Urban3. In this image we have selected only a region,
region A, with corners and edges. 
In Figure \ref{Grove2Urban3}(c) we plot the results for this region, when
$\lambda=0.9$ and $\theta=0.7$.
The best $\lambda_{SB}$ value is in this case $1$. We remind  the value of 
 $\lambda_{SB}$ seems to be related to the value of $\theta$, that is,  it should be close to 
 $2/\theta$.
The best results are attained  for $\alpha$ between $1.5$ and $2$ and in particular
for  $\alpha=1.5$.
The variations of the errors in terms of $\alpha$ are not smooth as in the case of the sequence Grove2
for the region A shown in Figure \ref{Grove2Urban3}(a).

\begin{figure}[h]
\includegraphics[height=40mm,width=56mm]
{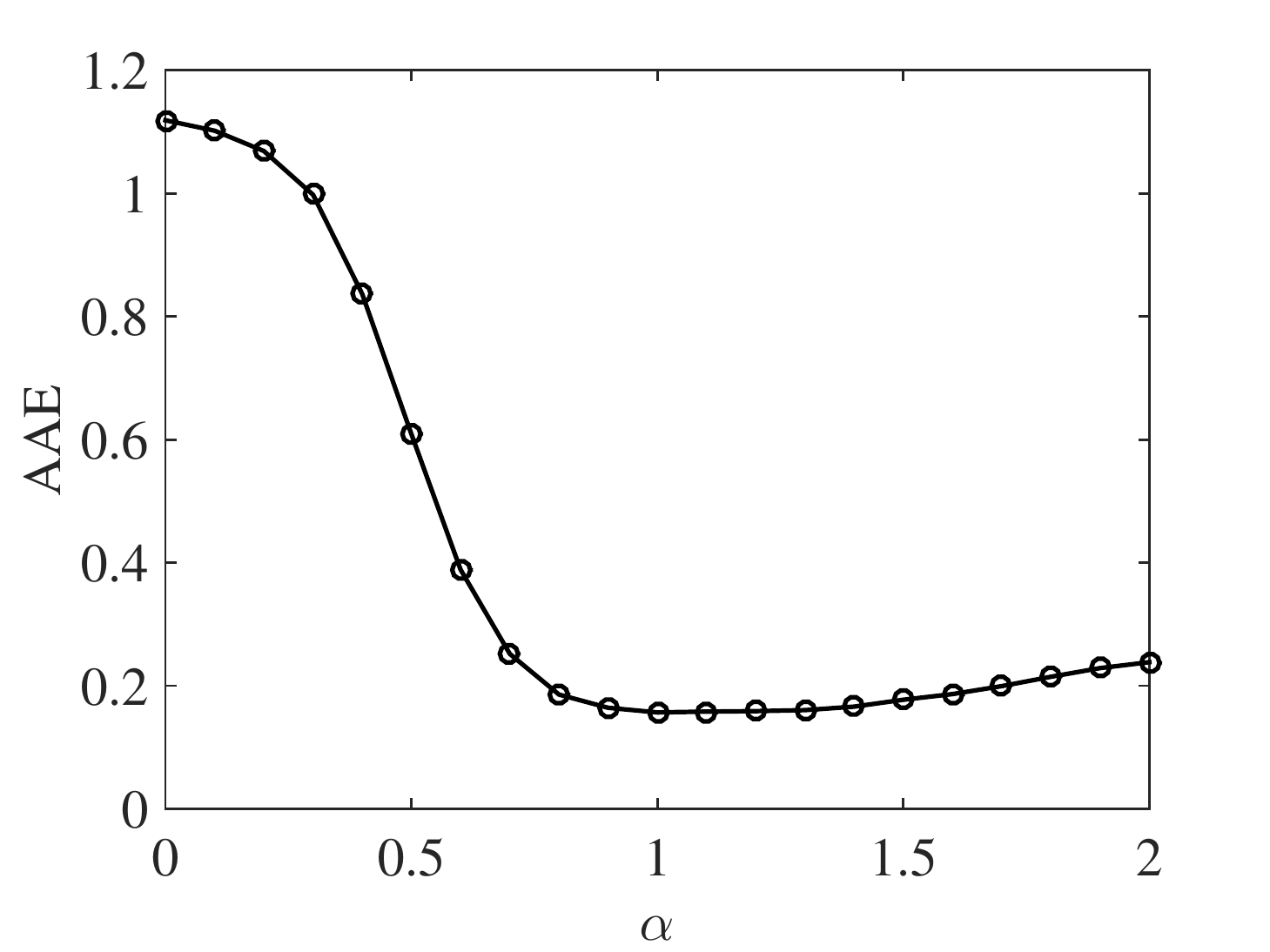}
\includegraphics[height=40mm,width=56mm]
{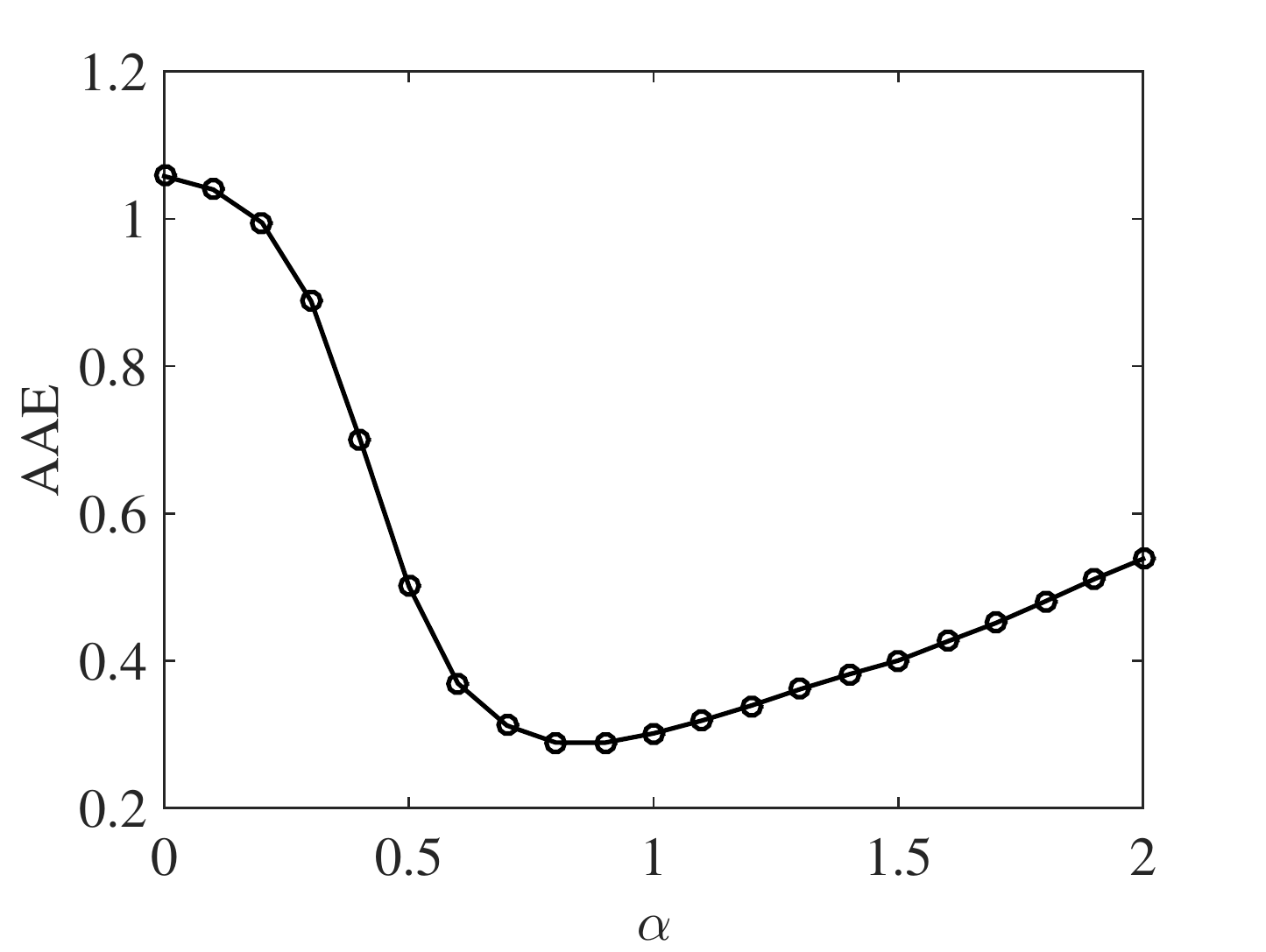}
\includegraphics[height=40mm,width=56mm]
{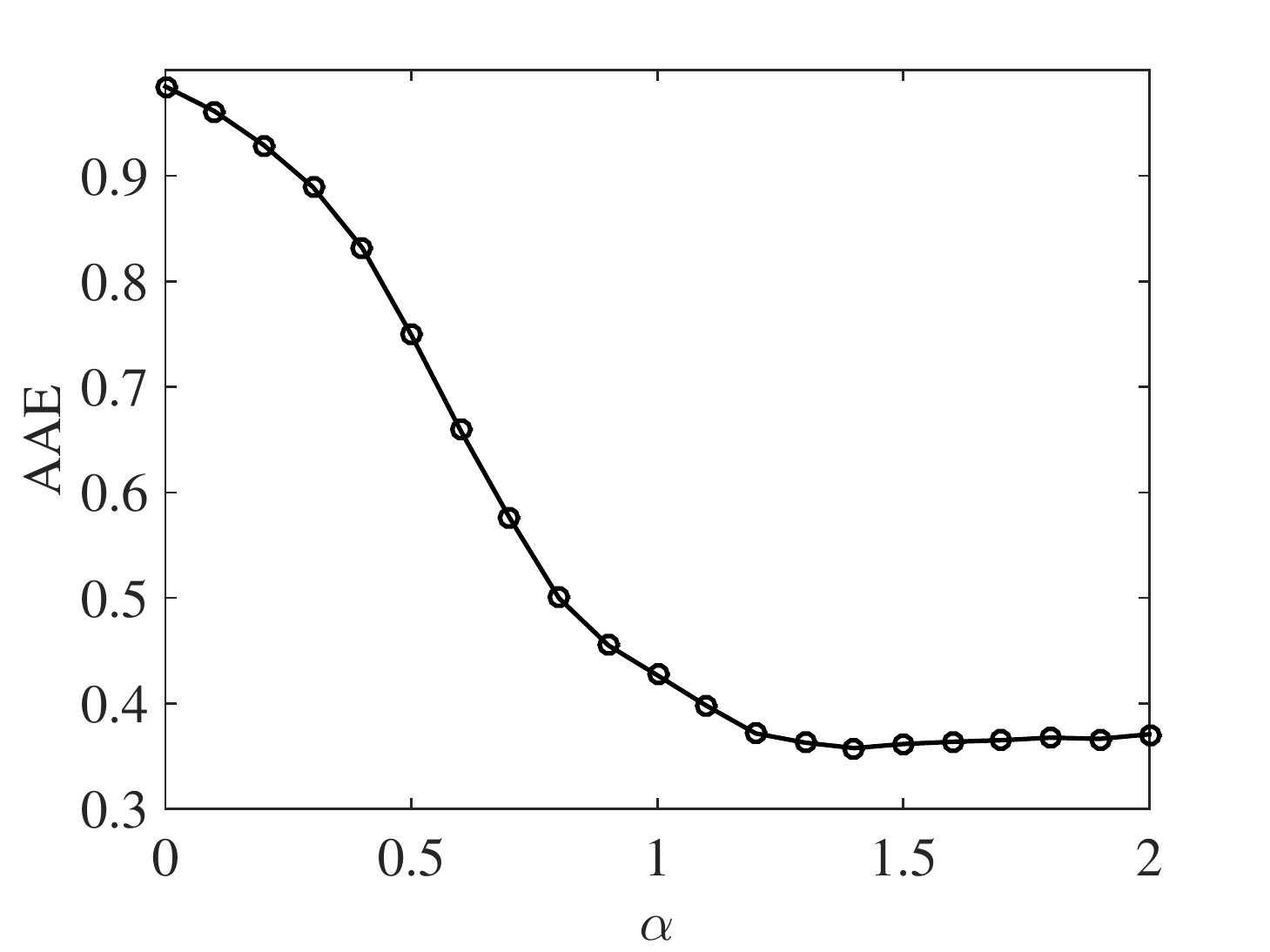}
\newline
\includegraphics[height=40mm,width=56mm]
{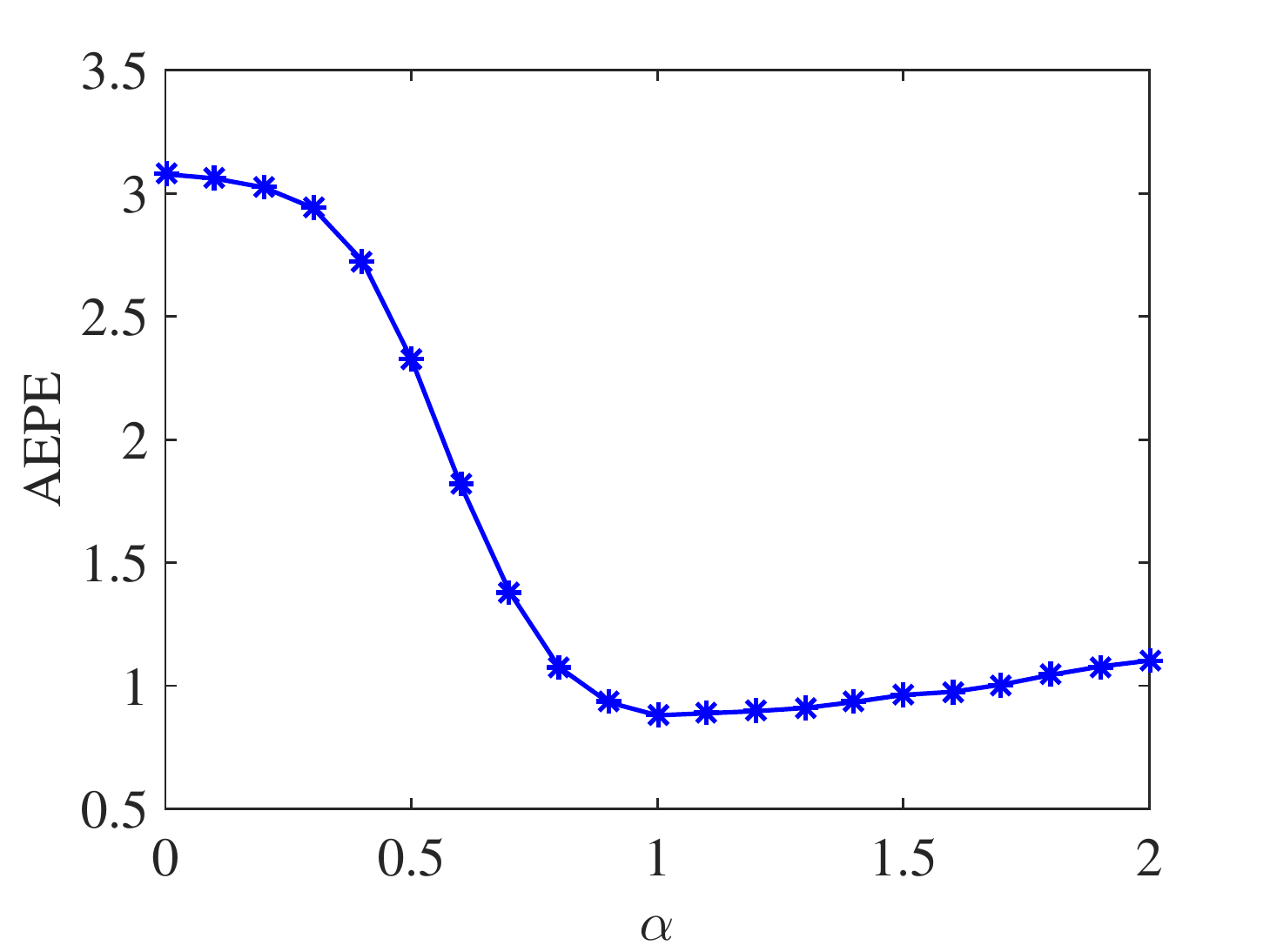}
\includegraphics[height=40mm,width=56mm]
{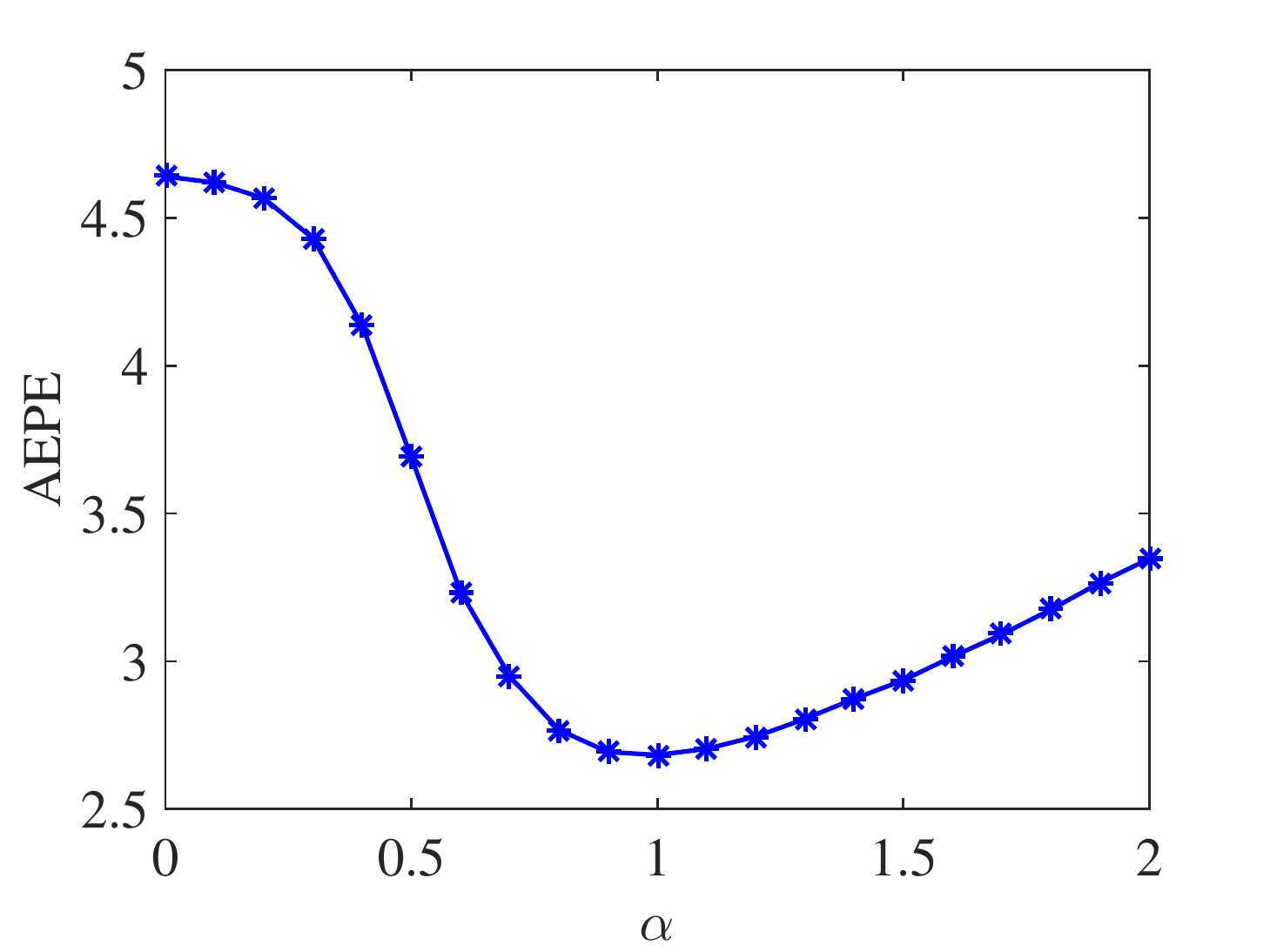}
\includegraphics[height=40mm,width=56mm]
{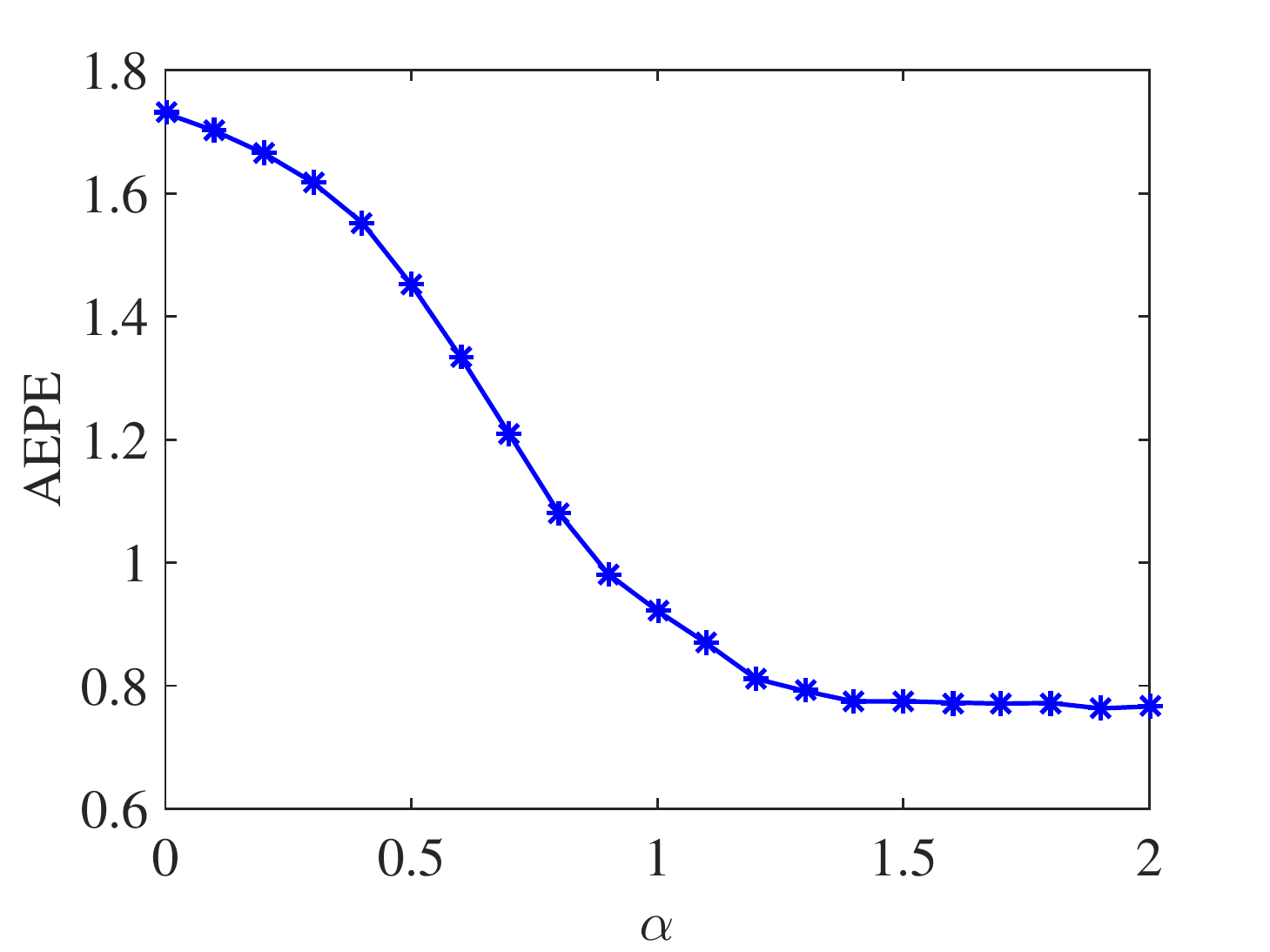}
{\small \hspace*{2.6cm} (a) \hspace*{5.2cm}(b)  \hspace*{5.2cm}(c)}
\caption{(a) Hydrangea -- PartA. For  $\lambda= 0.1$ and $\theta = 0.8$.
Best results for $\lambda_{SB}=3$. AAE: best result
  for $\alpha=1$. AEPE: best result for $\alpha=1$;
 (b) Hydrangea -- PartB. For  $\lambda= 0.1$ and $\theta = 0.8$.
Best results for $\lambda_{SB}=10$. AAE: best result
  for $\alpha=0.9$. AEPE: best result for $\alpha=1$;
 (c) Rubber Whale -- PartA. For  $\lambda= 0.4$ and $\theta = 0.4$.
Best results for $\lambda_{SB}=2$. AAE: best result
 for $\alpha=1.4$. AEPE: best result for $\alpha=1.4$.}
\label{HydraRubber}   
\end{figure}

In Figure \ref{HydraRubber} we present the results for  two other sequences:  
Hydrangea, Figures \ref{HydraRubber}(a) and  \ref{HydraRubber}(b),
and Rubber Whale, Figure \ref{HydraRubber}(c).  The regions analysed are marked  
in Figure \ref{fracimage}.  For  Hydrangea,  regions A and B, the errors are smaller
for values of $\alpha$ around $1$ and for  Rubber Whale, region A,
the smaller errors are for $1.4$. In the Rubber Whale case although the best value
is reached at $\alpha=1.4$ and not $\alpha=1$, 
we note that for these values of $\alpha$ the differences between the errors
are less relevant than for the sequences Grove2 and Urban3.

We have seen that for flat regions the best $\alpha$ is close to $0$, see Figure \ref{Grove2Urban3}(b),
for edges is between $1$ and $1.5$, see Figures \ref{Grove2Urban3}(a) and \ref{HydraRubber}(c),
and for corners between $1.5$ and $2$, see Figure \ref{Grove2Urban3}(c).
Additionally the results for the  most difficult datasets, Grove2 and Urban3,
highlights  the advantage of using the fractional order $\alpha$, by presenting
significantly smaller errors  for values of  $\alpha$ that are different from $1$.

\

\section{Final Remarks}

We have presented and tested  a new algorithm to solve an optical flow model. The novelty
of this work is twofold: there is the inclusion of a regularisation operator
which uses fractional derivatives and the application 
of  the split Bregman technique.

It is difficult to recover accurately motion fields and these difficulties arise from scene geometry
and texture complexity. 
We have seen that the parameter $\alpha$, related to the order of the fractional
regularisation operator, can be adjusted to deal with different regions. This motivates
a future work, to further explore the potential of the proposed approach,
that is  to build a feature based algorithm.
This algorithm would present a robust approach
that integrates region tracking, that is, it
would be able to detect {\it a priori} the type of structures the image  presents and
to adjust the parameter $\alpha$ accordingly.
Additionally, faster iterative methods can be developed to solve the iterative system
related to the fractional Euler-Lagrange equation.


\end{document}